\documentclass[aop,citesort,MSNbibl,dvips]{arximspdf}
\usepackage{graphics}
\doi{10.1214/09-AOP518}
\volume{38}
\issue{4}
\pubyear{2010}
\firstpage{1609}
\lastpage{1638}

\makeatletter

\newproclaim{Definition}{Definition}[section]
\newtheorem{Proposition}[Definition]{Proposition}
\newtheorem{Theorem}[Definition]{Theorem}
\newtheorem{Corollary}[Definition]{Corollary}
\newtheorem{Lemma}[Definition]{Lemma}
\newproclaim{Remark}[Definition]{Remark}

\makeatother

\begin{document}
\begin{frontmatter}

\title{Gibbsianness and non-Gibbsianness in divide and~color~models}
\runtitle{Gibbsianness and non-Gibbsianness in $\operatorname{DaC}(q)$ models}

\begin{aug}
\author[A]{\fnms{Andr\'as} \snm{B\'alint}\corref{}\ead[label=e1]{abalint@few.vu.nl}}
\runauthor{A. B\'alint}
\affiliation{VU University Amsterdam}
\address[A]{Department of Mathematics\\
VU University Amsterdam\\
De Boelelaan 1081a\\
1081 HV Amsterdam\\
The Netherlands\\
\printead{e1}} %adresu isvedimo komanda gale!
\end{aug}

\pdfauthor{Andras Balint}

% HISTORY:
\received{\smonth{4} \syear{2009}}
\revised{\smonth{11} \syear{2009}}

% ABSTRACT
%
\begin{abstract}
For parameters $p\in[0,1]$ and $q>0$ such that the Fortuin--Kasteleyn (FK) random-cluster
measure $\Phi_{p,q}^{\mathbb{Z}^d}$ for $\mathbb{Z}^d$ with parameters
$p$ and $q$ is unique, the $q$-divide and color [$\operatorname{DaC}(q)$] model on
$\mathbb{Z}^d$ is defined as follows. First, we draw a bond
configuration with distribution $\Phi_{p,q}^{\mathbb{Z}^d}$. Then, to
each (FK) cluster (i.e., to every vertex in the FK
cluster), independently for different FK clusters, we assign a spin
value from the set $\{1,2,\ldots,s\}$ in such a way that spin $i$ has
probability $a_i$.

In this paper, we prove that the resulting measure on spin
configurations is a Gibbs measure for small values of $p$ and is not a
Gibbs measure for large $p$, except in the special case of $q\in
\{2,3,\ldots\}$, $a_1=a_2=\cdots=a_s=1/q$, when the $\operatorname{DaC}(q)$ model
coincides with the $q$-state Potts model.
\end{abstract}

% KEYWORDS
%

\end{frontmatter}

%s1 ###
\section{Introduction}\label{intro}

The random-cluster representations of various models have played an
important role in the study of physical systems and phase transitions.
They provide a different viewpoint for physical models and many
problems in the Ising and Potts models can indeed be solved
by using their random-cluster representations
(see, e.g., \cite{Georgii,RCrepresentations,GHM}).

For $\beta\geq0$ and an integer $q\geq2$, a spin configuration in the
$q$-state Potts model at inverse temperature $\beta$
can be obtained as follows.
Draw a bond configuration according to a random-cluster measure with parameters
$p=1-e^{-2\beta}$ and $q$ (for definitions, see Section \ref
{definitions}), then assign to each vertex a spin
value from the set $\{1,2,\ldots,q\}$
in such a way that all spins have equal probability and that vertices
that are connected in the bond configuration get the same spin.
If the spin is chosen from a set $\{1,2,\ldots,s\}$ with an integer
$1<s<q$ and the probability of spin $i$ is $k_i/q$ with positive
integers $k_1,k_2,\ldots,k_s$ such that $\sum_{i=1}^sk_i=q$, then we
get the so-called fuzzy Potts model \cite{MVV,fuzzyGibbs}.
Recent papers (see \cite{fuzzyposcorr,DaC,Garet,KW,BCM,Ising,GrGr})
have shown that generalizations of the above constructions
with different values of $q$ and $s$, as well as more general rules of
spin assignment,
are also of interest. From a mathematical viewpoint,
such models are natural examples of a dependent site percolation model
with a simple definition but nontrivial behavior.
The study of such models may also lead to a better understanding of
models of primary physical importance, as was the case in \cite{Ising},
where an informative, new perspective on the high temperature Ising
model on the triangular lattice was given.

The model treated here is defined as follows. Let $G=(\mathcal
{V},\mathcal{E})$ be a (finite or infinite) locally finite graph.
Fix parameters $p\in[0,1]$, $q>0$ in such a way that
there exists exactly one random-cluster measure for $G$ with parameters
$p$ and $q$.
We denote this measure by $\Phi_{p,q}^G$.
Also, fix an integer $s\geq2$ and $a_1,a_2,\ldots,a_s\in(0,1)$ such
that $\sum_{i=1}^s a_i=1$,
and define the \textit{single-spin space} $S=\{1,2,\ldots,s\}$ and the
\textit{state space} $\Omega^G=\Omega_C^G\times\Omega_D^G$ with
$\Omega_C^G=S^{\mathcal{V}}$ and $\Omega_D^G=\{0,1\}^{\mathcal{E}}$.
Let $Y$ be a random bond configuration taking values in $\Omega_D^G$
with distribution $\Phi_{p,q}^{G}$.
Given $Y=\eta$ for some $\eta\in\Omega_D^G$, we construct a random
$\Omega_C^G$-valued spin configuration $X$ by assigning spin $i\in S$
with probability $a_i$ to each connected component in $\eta$ (i.e.,
the same spin $i$ to each vertex in the component), independently for
different components.
We write $\mathbb{P}^{G}_{p,q,(a_1,a_2,\ldots,a_s)}$ for the
joint distribution of $(X,Y)$ on $\Omega^G$
and $\mu^{G}_{p,q,(a_1,a_2,\ldots,a_s)}$ for the marginal\vspace*{1pt} of $\mathbb
{P}^{G}_{p,q,(a_1,a_2,\ldots,a_s)}$ on $\Omega_C^G$.
This definition is a slight generalization of the fractional fuzzy
Potts model defined in \cite{fuzzyposcorr}, page~1156 (see also \cite
{BCM}, Section 1.2). However, we shall call this model the
\textit{$q$-divide and color} [$\operatorname{DaC}(q)$] \textit{model} to emphasize that we look
at it as a generalization of the model introduced in \cite{DaC}
by H\"aggstr\"om [which is the DaC(1) model in the present
terminology], rather than
of the fuzzy Potts model of \cite{MVV,fuzzyGibbs}.

Let us now consider the (hypercubic) lattice with vertex set $\mathbb
{Z}^d$ and edge set $\mathcal{E}^d$ with edges between vertices at
Euclidean distance $1$. With an abuse of notation, we shall denote this
graph by $\mathbb{Z}^d$ and the sets $\Omega^{\mathbb{Z}^d}_D,\Omega
^{\mathbb{Z}^d}_C$ and $\Omega^{\mathbb{Z}^d}$ by $\Omega_D,\Omega
_C$ and $\Omega$, respectively.
The present work is focused on the Gibbs properties and
$k$-Markovianness of the measure $\mu^{\mathbb
{Z}^d}_{p,q,(a_1,a_2,\ldots,a_s)}$ in $d\geq2$ dimensions.\vspace*{1pt} Since the
cases $p=0$ and $p=1$ are trivial, we henceforth assume that $p\in(0,1)$.
We only give results for $q\geq1$
since much more is known about random-cluster measures with $q\geq1$
than with $q<1$.

We shall prove that, except in the special case of $q=s$ and
$a_1=a_2=\cdots=a_s$ [when the $\operatorname{DaC}(q)$ model coincides with the
$q$-state Potts model on $\mathbb{Z}^d$ at inverse temperature $\beta
=-1/2\log(1-p)$], the $\operatorname{DaC}(q)$ model is not $k$-Markovian for any~$k$.
For large values of $p$, $\mu^{\mathbb{Z}^d}_{p,q,(a_1,a_2,\ldots
,a_s)}$ is not even\vspace*{2pt} quasilocal and is therefore not a Gibbs measure,
again with the exception of the Potts case.
This shows that the Gibbsianness of the Potts model at low temperatures
is very sensitive with respect to perturbations in the assignment of
the spin probabilities; even the smallest change makes the model nonquasilocal.
By demonstrating the special role of the $q$-state Potts Gibbs measure
among $\operatorname{DaC}(q)$ models, our result supports the view expressed in \cite
{vEFS2,vEFS} that, especially at low temperatures, Gibbsianness of
measures is the exception rather than the rule. For a related general
result, see \cite{Israel}, where Israel proved that in the set of all
translation invariant measures,
Gibbsianness is exceptional in a topological sense.
However, if $p$ is small enough,
then Gibbsianness does hold.
The proof of this fact uses the idea
that at small values of $p$ [which correspond to high temperatures in
the (fuzzy) Potts model], the $\operatorname{DaC}(q)$ model
is close in spirit to independent site percolation on $\mathbb{Z}^d$.

These results are in line with
those in \cite{fuzzyGibbs} (see also \cite{MVV} and \cite{HK})
concerning the fuzzy Potts model and, in some cases, essentially the
same proofs work in the current, more general situation. Therefore, in
some instances, only a sketch of the proof is given and the reader is
referred to \cite{fuzzyGibbs} for the details.
Note, however, that such similarities are not immediate from the
definitions of the models.
More importantly, in the $\operatorname{DaC}(q)$ model, a distinction must be made
between the case when $a_i\geq1/q$ for all $i$ and when there exists
some $j$ with $a_j<1/q$. In the former (of which the fuzzy Potts model
is a special case), a rather complete picture can be given, whereas in
the latter, there is an interval in $p$ where we do not know whether
$\mu^{\mathbb{Z}^d}_{p,q,(a_1,a_2,\ldots,a_s)}$ is a Gibbs
measure.\vspace*{2pt}

Finally, we give a sufficient (but not necessary) condition for the
almost sure quasilocality of $\mu^{\mathbb{Z}^d}_{p,q,(a_1,a_2,\ldots
,a_s)}$ and,\vspace*{1pt} as an application, we obtain this weak form of
Gibbsianness in the two-dimensional case for a large range of
parameters. Some intuition underlying our main results will be
described after Remark \ref{PottsMarkovian}.

%s2 ###
\section{Definitions and main results}\label{definitions}

%s2.1 ###
\subsection{Random-cluster measures}\label{rcmeasures}

In this section, we recall the definition of Fortuin--Kasteleyn (FK)
random-cluster measures and those properties of these measures that
will be important throughout the rest of the paper. For the proofs and
much more on random-cluster measures, see, for example, \cite{grimmett2}.
\begin{Definition}\label{FKdeffin}
For a finite graph $G=(\mathcal{V},\mathcal{E})$ and parameters $p\in
[0,1]$ and $q>0$, the \textit{random-cluster measure} $\Phi^G _{p,q}$ is
the measure on $\Omega_D^G$
which assigns to a bond configuration $\eta\in\Omega_D^G$
the probability
%
%e1 ###
\begin{equation}\label{FKfinite}
\Phi^G _{p,q}(\eta)= \frac{q^{k(\eta)}}{Z^G_{p,q}}\prod_{e\in
\mathcal{E}}p^{\eta(e)}(1-p)^{1-\eta(e)},
\end{equation}
where $k(\eta)$ is the number of connected components
in the graph with vertex set $\mathcal{V}$ and edge set $\{e\in\mathcal
{E}\dvtx \eta(e)=1\}$
(we call such components \textit{FK clusters} throughout, edges with
state $1$ \textit{open}
and edges with state $0$ \textit{closed}),
and $Z^G_{p,q}$ is the appropriate normalizing factor.
\end{Definition}

This definition is not suitable for infinite graphs.
In that case, we shall
require that certain conditional probabilities are the same as in the
finite case.
The relevant definition, given below, will formally contain
conditioning on an event with probability $0$, which should be
understood as conditioning on the appropriate $\sigma$-algebra. We
shall frequently use this simplification in order to keep the notation
as simple as possible.

A graph is called \textit{locally finite} if every vertex has a bounded degree.
We shall denote bond configurations throughout by $\eta$ and $\zeta$.
For the restriction of a bond configuration $\eta$ to an edge set $H$,
we write $\eta_H$.
For vertices $v$ and $w$, we denote the edge between $v$ and $w$ by
$ \langle v,w \rangle  $.
The following definition is taken from \cite{fuzzyGibbs} and its
equivalence with a more common definition (where arbitrary finite edge
sets and not only single edges are considered) is stated, for example,
in Lemma 6.18 of~\cite{GHM}.
\begin{Definition}\label{FKdef}
For
an infinite, locally finite graph $G=(\mathcal{V},\mathcal{E})$ and
parameters $p\in[0,1]$, $q>0$,
a measure $\phi$ on $\Omega_D^G$
is called a \textit{random-cluster measure for $G$ with parameters $p$
and $q$} if,
for each edge $e= \langle x,y \rangle \in\mathcal{E}$ and edge
configuration $\zeta\in\{0,1\}^{\mathcal{E}\setminus\{e\}}$ outside
$e$, we have that
\[
\phi\bigl(\{\eta\in\Omega_D^G\dvtx\eta(e)=1\}\mid\bigl\{\eta\in\Omega
_D^G\dvtx\eta_{\mathcal{E}\setminus\{e\}}=\zeta\bigr\}\bigr)=
\cases{
p, &\quad if $\stackrel{\zeta}{x\leftrightarrow y}$,\cr
\dfrac{p}{p+(1-p)q}, &\quad otherwise,}
\]
where $\stackrel{\zeta}{x\leftrightarrow y}$ denotes that there exists
a path of edges between $x$ and $y$ in which every edge has $\zeta
$-value $1$.
\end{Definition}

It is not difficult to prove that one gets the same
conditional probabilities for random-cluster measures on finite graphs,
so Definition \ref{FKdef} is a reasonable extension of Definition \ref
{FKdeffin} to infinite graphs.

It is not clear from the definition that such measures exist. However,
for $\mathbb{Z}^d$ and $q\geq1$,
two random-cluster measures can be constructed as follows.
For a vertex set $H\subset\mathbb{Z}^d$, let $\partial H$ denote the
\textit{vertex boundary} of the set, that is, $\partial H=\{v\in\mathbb
{Z}^d\setminus H\dvtx\exists w\in H$ such that $ \langle v,w \rangle  \in\mathcal{E}^d\}$.
Define, for $n\in\{1,2,\ldots\}$, the set $\Lambda_n=\{-n,\ldots,n\}
^d$ and the graph $G_n=(\mathcal{V}_n,\mathcal{E}_n)$ with vertex set
$\mathcal{V}_n=\Lambda_n\cup\partial\Lambda_n$
and edge set $\mathcal{E}_n=\{e\in\mathcal{E}^d$: both endvertices of
$e$ are in $\mathcal{V}_n\}$. For $n\in\{1,2,\ldots\}$, let $W_n$ be
the event that all edges with both endvertices in $\partial\Lambda_n$
are open
and let $\Phi^{G_n,1} _{p,q}$ be the measure $\Phi^{G_n} _{p,q}$
conditioned on $W_n$.
Then both $\Phi^{G_n} _{p,q}$ and $\Phi^{G_n,1} _{p,q}$ converge
weakly as $n\to\infty$; we denote the limiting measures by $\Phi
_{p,q}^{\mathbb{Z}^d,0}$ and $\Phi_{p,q}^{\mathbb{Z}^d,1}$,
respectively. $\Phi_{p,q}^{\mathbb{Z}^d,0}$ is called the \textit{free},
and $\Phi_{p,q}^{\mathbb{Z}^d,1}$ the \textit{wired}, random-cluster
measure for $\mathbb{Z}^d$ with parameters $p$ and $q$.
These measures are indeed random-cluster measures in the sense of
Definition \ref{FKdef}, moreover, they are extremal among such measures
in the following sense.

A natural partial order on the set $\Omega_D=\{0,1\}^{\mathcal{E}^d}$ of
edge configurations is given by defining $\eta^{\prime}\geq\eta$
for $\eta,\eta^{\prime}\in\Omega_D$
if, for all $e\in\mathcal{E}^d$, $\eta^{\prime}(e)\geq\eta(e)$.
We call a function $f\dvtx\Omega_D\to\mathbb{R}$
\textit{increasing} if $\eta^{\prime}\geq\eta$ implies that $f(\eta
^{\prime})\geq f(\eta)$.
For probability measures $\phi,\phi^{\prime}$ on $\Omega_D$,
we say that $\phi^{\prime}$ is \textit{stochastically larger} than $\phi
$ if, for all bounded increasing measurable
functions $f\dvtx\Omega_D\to\mathbb{R}$,
we have that
\[
\int_{\Omega_D}f(\eta) \,d\phi^{\prime}(\eta)\geq\int_{\Omega
_D}f(\eta) \,d\phi(\eta).
\]
For later purposes, we remark that
by Strassen's theorem \cite{Strassen},
this is equivalent to the existence of an appropriate coupling of the
measures $\phi^{\prime}$ and $\phi$,
that is, the existence of a probability measure $Q$ on $\Omega_D\times
\Omega_D$ such that the marginals of $Q$ on the first and second
coordinates are $\phi^{\prime}$ and $\phi$, respectively, and $Q(\{
(\eta^{\prime},\eta)\in\Omega_D\times\Omega_D\dvtx\eta^{\prime
}\geq\eta\})=1$.

It is well known that
$\Phi_{p,q}^{\mathbb{Z}^d,0}$ is the stochastically smallest, and
$\Phi_{p,q}^{\mathbb{Z}^d,1}$ the stochastically largest,
random-cluster measure for $\mathbb{Z}^d$ with parameters $p$ and $q$.
Therefore, there exists a unique random-cluster measure for $\mathbb
{Z}^d$ with parameters $p$ and $q$ if and only if
%
%e2 ###
\begin{equation}\label{uniqueFK}
\Phi_{p,q}^{\mathbb{Z}^d,0}=\Phi_{p,q}^{\mathbb{Z}^d,1}.
\end{equation}
This is the case for any fixed $q\geq1$, except (possibly) for at most
countably many values of $p$.
It is widely believed that for any $q\geq1$, there is at most one
exceptional~$p$, which can only be
the \textit{critical value}
$p_c(q,d)=\sup\{p\dvtx\Phi_{p,q}^{\mathbb{Z}^d,0}(\{\eta\in\Omega
_D\dvtx\mathbf{0}$ is in an infinite FK cluster in $\eta\})=0\}$,
where \textbf{0} denotes the origin in $\mathbb{Z}^d$.
It is not difficult to show that the choice of $\Phi_{p,q}^{\mathbb
{Z}^d,0}$ in the definition is not crucial. That is, for
any random-cluster measure $\phi$ for $\mathbb{Z}^d$ with parameters
$p$ and $q$, we have that
\[
\phi(\{\eta\in\Omega_D\dvtx\mathbf{0}\mbox{ is in an infinite FK cluster
in }\eta\} )
\cases{
=0, &\quad if $p<p_c(q,d)$,\cr
>0, &\quad if $p>p_c(q,d)$.}
\]
For the rest of the paper, we will assume, without further mention,
that the parameters $d,p,q$
for the $\operatorname{DaC}(q)$ model on $\mathbb{Z}^d$ are always chosen in such a
way that (\ref{uniqueFK}) holds,
and we will denote the unique random-cluster measure by $\Phi^{\mathbb
{Z}^d}_{p,q}$.

Another important feature of
the random-cluster measures $\Phi^{\mathbb{Z}^d,0}_{p,q}$ and $\Phi
^{\mathbb{Z}^d,1}_{p,q}$ with $q\geq1$ is that they satisfy the FKG inequality
for increasing events \cite{FKG} (an event $A\subset\Omega_D$ is called
\textit{increasing} if its indicator
function is increasing, that is, if $\eta\in A$ and $\eta^{\prime
}\geq\eta$ implies that $\eta^{\prime}\in A$).
This, in particular, means that for
$d\geq2$, $p\in[0,1]$, $q\geq1$,
any edge set $E\subset\mathcal{E}^d$,
configuration $\zeta\in\{0,1\}^E$ on $E$ and
increasing events $A_1,A_2\subset\Omega_D$, we have, letting
$B=\{ \eta\in\Omega_D\dvtx\eta_{E}=\zeta\}$, that
%
%e3 ###
\begin{equation}\label{strongFKGfree}
\Phi^{\mathbb{Z}^d,0}_{p,q} (A_1\cap A_2\mid B) \geq\Phi^{\mathbb
{Z}^d,0}_{p,q}(A_1\mid B)\Phi^{\mathbb{Z}^d,0}_{p,q}(A_2\mid B).
\end{equation}

Finally, for our main result, we also need to consider the critical
value in half-spaces.
Let $\mathcal{H}^+=\mathcal{H}^+_d$ denote the subset of $\mathbb{Z}^d$
which consists of those vertices whose first coordinate is strictly
positive, let
$\tilde{E}\subset\mathcal{E}^d$
denote the set of edges
that are incident to at least one vertex in $\mathbb{Z}^d\setminus
\mathcal{H}^+$
and denote the vertex $(1,0,0,\ldots,0)\in\mathbb{Z}^d$ by $u_1$.
Also, consider the event
$A_{\mathcal{H}^+}=\{\eta\in\Omega_D\dvtx u_1$ is in an infinite open
path in $\eta$ which is contained in $\mathcal{H}^+\}$.
For $q\geq1$, we define
$p_c^{\mathcal{H}}(q,d)=\sup\{p\dvtx\Phi_{p,q}^{\mathbb
{Z}^d,0}(A_{\mathcal{H}^+}\mid\{\eta\in\Omega_D\dvtx \eta_{\tilde
{E}}\equiv0\})=0\}$.

Using (\ref{strongFKGfree}), it is easy to see that $p_c^{\mathcal
{H}}(q,d)\geq p_c(q,d)$. Equality of the two critical values for $q=1$
was proven by Barsky, Grimmett and Newman \cite{BGN}, for $q=2$ by
Bodineau \cite{Bodineau} and for very large values of $q$, it follows from
the Pirogov--Sinai theory (see the last paragraph of Section 2.3 in \cite
{Bodineau}).
For general $q\geq1$, equality has been conjectured
\cite{Pisztora,fuzzyGibbs,Bodineau,Wouts}, but no
definite result has been established thus far. However, an upper bound
$
p_c^{\mathcal{H}}(q,d)\leq\frac{p_c(1,d)q}{p_c(1,d)q+1-p_c(1,d)}
$
can easily be given, using the fact that for $q\geq1$, $\Phi
_{p,q}^{\mathbb{Z}^d,0}$ conditioned on $\{\eta\in\Omega_D\dvtx \eta
_{\tilde{E}}\equiv0\}$ is stochastically larger on $\mathcal
{E}^d\setminus\tilde{E}$ than
$\Phi_{{p}/({p+(1-p)q}),1}^{\mathbb{Z}^d,0}$,
and the fact that $p_c^{\mathcal{H}}(1,d)=p_c(1,d)$.
Note that $p_c(1,d)$
is the critical value for Bernoulli bond percolation on $\mathbb{Z}^d$.
It is well known (see, e.g., \cite{Grimmett}) that for all $d\geq2$,
$0<p_c(1,d)<1$. This implies that the above upper bound for
$p_c^{\mathcal{H}}(q,d)$ is nontrivial.

%s2.2 ###
\subsection{Main results}\label{mainresults}

Before stating the main results, let us give the relevant definitions.
In this section, $\mu$ denotes a probability measure on
$\Omega_C={S}^{\mathbb{Z}^d}$. Spin configurations will be denoted
throughout by $\xi,\sigma$ and $\kappa$, and
the restriction of a spin configuration $\xi$ to a vertex set $W$ by
$\xi_W$.
For a set $W\subset\mathbb{Z}^d$ and a spin configuration $\sigma\in
S^W$ on $W$, we define $K_W ^{\sigma}=\{\xi\in\Omega_C\dvtx \xi
_W=\sigma\}$.
We shall use $A\subset\subset B$ to denote that ``$A$ \textit{is a finite
subset of} $B$'' throughout. We denote the graph theoretic distance on
$\mathbb{Z}^d$ by $\mathrm{dist}$ and define the \textit{distance} between a vertex
$v\in\mathbb{Z}^d$ to a vertex set $H\subset\mathbb{Z}^d$ by
$\operatorname{dist}(v,H)=\min\{\operatorname{dist}(v,w)\dvtx w\in H\}$.
For
$k\in\{1,2,\ldots\}$, let $\partial_kH$ denote the
\textit{$k$-neighborhood} of $H$, that is, $\partial_kH=\{v\in\mathbb
{Z}^d\dvtx 1\leq \operatorname{dist}(H,v)\leq k\}$. Note that $\partial_1H=\partial H$.

We usually want to view the $\operatorname{DaC}(q)$ model as a dependent spin model on
$\mathbb{Z}^d$, in which
the only role of the edge configuration is to introduce the dependence.
One of the first questions which naturally arises concerning a spin
model is whether the finite energy property of \cite{NSch} holds. This
turns out to be the case; moreover, we can even prove a stronger form
of it, called uniform nonnullness. The proofs of all statements in
this section will be given in Section \ref{mainproofs}.
\begin{Definition}\label{undef}
$\mu$ is called \textit{uniformly nonnull} if there exists an
$\varepsilon>0$ such that for all $v\in\mathbb{Z}^d$, $m\in{S}$ and
$\sigma\in{S}^{\mathbb{Z}^d\setminus\{v\}}$, we have that
\[
\mu\bigl(K_{\{v\}}^m\mid K_{\mathbb{Z}^d\setminus\{ v\}}^{\sigma}\bigr)\geq
\varepsilon.
\]
\end{Definition}
\begin{Proposition}\label{untheorem}
For all $d\in\{1,2,\ldots\}$, $q\geq1$, $p\in[0,1)$ and arbitrary
values of the other parameters,
the measure $\mu^{\mathbb{Z}^d}_{p,q,(a_1,\ldots,a_s)}$ is uniformly nonnull.
\end{Proposition}

The concept of $k$-Markovianness is concerned with the following question:
conditioning on a spin configuration outside a set $W$, do
vertices farther than $k$ from $W$ have any influence on
the spin configuration in $W$?
\begin{Definition}
For $k\in\{1,2,\ldots\}$, $\mu$ is called \textit{$k$-Markovian} if,
for all $W\subset\subset\mathbb{Z}^d$, $\kappa\in{S}^{W}$ and
$\sigma,\sigma^{\prime}\in{S}^{\mathbb{Z}^{d}\setminus W}$ such
that $\sigma_{\partial_kW}=\sigma^{\prime}_{\partial_kW}$, we have that
\[
\mu(K _W^{\kappa}\mid K_{\mathbb{Z}^d\setminus W}^{\sigma})=\mu
(K_W^{\kappa}\mid K_{\mathbb{Z}^d\setminus W}^{\sigma^{\prime}} ).
\]
\end{Definition}

A weaker notion is that of quasilocality, where the above conditional
probabilities do not need to be equal for any $k$; the only requirement
is that their difference tends to $0$ as $k\to\infty$. Due to the
compactness of ${S}^{\mathbb{Z}^d}$ in the product topology, this
amounts to the following definition.
\begin{Definition}\label{quasilocality}
$\mu$ is called \textit{quasilocal} if, for all $W\subset\subset\mathbb
{Z}^d$, $\kappa\in{S}^W$ and $\sigma\in{S}^{\mathbb{Z}^{d}\setminus
W}$, we have that
\[
\lim_{k\to\infty}\mathop{\sup_{\sigma^{\prime}\in\mathbb
{Z}^d\setminus W}}_{\sigma^{\prime} _{\partial_kW}=\sigma_{\partial
_kW}}|\mu(K_W^{\kappa}\mid K_{\mathbb{Z}^d\setminus W}^{\sigma})-\mu
(K_W^{\kappa}\mid K_{\mathbb{Z}^d\setminus W}^{\sigma^{\prime}} )|=0.
\]
If the above equation holds for $\mu$-almost all $\sigma\in
{S}^{\mathbb{Z}^{d}\setminus W}$, then $\mu$ is called \textit{almost
surely quasilocal}.
\end{Definition}

Finally, we need to say what we mean by Gibbsianness.
Instead of the usual definition with absolutely summable interaction
potentials (see, e.g., \cite{Georgii,vEFS}), we shall use a well-known
characterization (see \cite{vEFS}, Theorem 2.12), namely that $\mu$ is
a \textit{Gibbs measure} if and only if it is quasilocal and uniformly nonnull.

We are now ready to state our main result concerning $k$-Markovianness
and Gibbsianness of the $\operatorname{DaC}(q)$ model. The cases $p=0,1$ are trivial,
so we assume that $p\in(0,1)$. For fixed $q$, $s$ and $a_1,\ldots
,a_s$, recall that $S=\{1,2,\ldots,s\}$ and define $S_{1/q}=\{i\in S\dvtx
a_i=1/q\}$.
The case $S=S_{1/q}$ is well understood since $S=S_{1/q}$ implies that
$s=q$ and $a_1=a_2=\cdots=a_s$,
in which case the procedure defining the $\operatorname{DaC}(q)$ model gives the
random-cluster representation of the Potts model.
Therefore, for all $p\in(0,1)$, $\mu^{\mathbb
{Z}^d}_{p,q,(1/q,1/q,\ldots,1/q)}$ equals a Gibbs measure
for the $q$-state Potts model on $\mathbb{Z}^d$ [at inverse temperature
$\beta=-1/2\log(1-p)$].
It follows immediately from the standard definition of Potts Gibbs measures
with a Hamiltonian (see, e.g., \cite{GHM} for the definition) that all
such measures are Markovian (i.e., \mbox{1-Markovian}).
For an alternative proof of the Markovianness of $\mu^{\mathbb
{Z}^d}_{p,q,(1/q,1/q,\ldots,1/q)}$,\vspace*{2pt} see Remark \ref{PottsMarkovian}.
If $S\neq S_{1/q}$, let $\ell\in S$ be an (for concreteness, the
smallest) index such that $a_{\ell}=\min\{a_i\dvtx i\in S\setminus S_{1/q}\}$.
\begin{Theorem}\label{mainthm}
Assume that $d\geq2$, $q\geq1$ and $S\neq S_{1/q}$. We then have the
following.
\begin{enumerate}
\item For any values of $p,a_1,\ldots,a_s\in(0,1)$, the measure $\mu
^{\mathbb{Z}^d}_{p,q,(a_1,\ldots,a_s)}$ is not $k$-Markovian for any
$k\in\{1,2,\ldots\}$.
\item If $a_{\ell}>1/q$, then:
\begin{enumerate}[(a)]
\item[(a)] for $p<p_c(qa_{\ell},d)$, $\mu^{\mathbb{Z}^d}_{p,q,(a_1,\ldots
,a_s)}$ is quasilocal; but,
\item[(b)] for
$p>p_c^{\mathcal{H}}(qa_{\ell},d)$, it is not quasilocal.
\end{enumerate}
\item If $a_{\ell}<1/q$, then:
\begin{enumerate}[(a)]
\item[(a)] if $p<\frac{p_c(1,d)qa_{\ell}}{p_c(1,d)qa_{\ell}+1-p_c(1,d)}$,
then $\mu^{\mathbb{Z}^d}_{p,q,(a_1,\ldots,a_s)}$ is quasilocal; whereas,
\item[(b)] if $p>p_c(1,d)$, it is not.
\end{enumerate}
\end{enumerate}
\end{Theorem}

Combining Theorem \ref{mainthm} with Theorem \ref{untheorem} and the
characterization of Gibbs measures mentioned earlier, we arrive at the
following result.
\begin{Corollary}\label{mainthmcor}
For $S=S_{1/q}$ and in cases 2\textup{(a)} and 3\textup{(a)} of Theorem \ref{mainthm},
the measure $\mu^{\mathbb{Z}^d}_{p,q,(a_1,\ldots,a_s)}$ is a
Gibbs\vspace*{2pt}
measure. However, in cases 2\textup{(b)} and 3\textup{(b)} of Theorem~\ref{mainthm},
$\mu _{p,q,(a_1, \ldots ,a_s)} ^{\mathbb{Z}^d}$ is not
a Gibbs measure.
\end{Corollary}

To demonstrate the fundamental difference between the $q$-state Potts
model and other $\operatorname{DaC}(q)$ models, let us consider the case with $s=q-1$
and $a_1=a_2=\cdots=a_{s}=\frac{1}{q-1}$.
Intuitively, for very large values of $q$, the difference between this
scenario and the case where $S=S_{1/q}$ should vanish. Nevertheless, while
$\mu^{\mathbb{Z}^d}_{p,q,(1/q,1/q,\ldots,1/q)}$\vspace*{2pt}
is a Gibbs measure for any $p$ and $q$, Corollary \ref{mainthmcor}
gives that
there exists a constant $c=c(d)\in(0,1)$ such that for all $q\in\{
3,4,\ldots\}$ and $p>c$,
$\mu^{\mathbb{Z}^d}_{p,q,({1}/({q-1}),{1}/({q-1}),\ldots,{1}/({q-1}))}$
is not\vspace*{2pt} a Gibbs measure.
This result might seem to contradict Theorem 2.9 in \cite{KO}, which
implies that any sufficiently fine local coarse graining preserves the
Gibbs property of the $q$-state Potts model.
Note, however, that an arbitrarily fine coarse graining is available
only when the local state space is continuous, which is not the case here.

The question of whether quasilocality is ``seriously'' violated in
cases where $\mu^{\mathbb{Z}^d}_{p,q,(a_1,\ldots,a_s)}$ is not
a\vspace*{1pt}
Gibbs measure
(i.e., whether ``bad'' configurations are exceptional or they actually
occur) is related to that of percolation
by the following statement, which is a generalization of Proposition
3.7 in \cite{DaC}.
\begin{Proposition}\label{almquas}
Consider the event $E_{\infty}=\{\xi\in\Omega_C\dvtx \xi$ contains an
infinite connected component of equal spins$\}$.
If the parameters $p\in[0,1]$, $q\geq1$,
$s\in\{2,3,\ldots\}$ and $a_1,\ldots,a_s\in(0,1)$ of the $\operatorname{DaC}(q)$
model are chosen in such a way that
%
%e4 ###
\begin{equation}\label{noperc}
\mu^{\mathbb{Z}^d}_{p,q,(a_1,\ldots,a_s)}(E_{\infty})=0,
\end{equation}
then $\mu^{\mathbb{Z}^d}_{p,q,(a_1,\ldots,a_s)}$
satisfies almost sure quasilocality.
\end{Proposition}

It is easy to see that (\ref{noperc}) is not a necessary condition for
almost sure quasilocality. For instance, one can take $d\geq2$, $q\geq
1$, $p=0$, $s=2$ and an $a_1<1$ which is greater than the critical
value for Bernoulli \textit{site} percolation on $\mathbb{Z}^d$. Then,
although (\ref{noperc}) fails,\vspace*{2pt} $\mu^{\mathbb{Z}^d}_{0,q,(a_1,a_2)}$ is
Markovian (and therefore obviously almost surely quasilocal). Despite
this, Proposition \ref{almquas} is not useless. We shall demonstrate
this below by giving an application in the two-dimensional case.
H\"aggstr\"om's results in Section 3 of \cite{fuzzyposcorr} imply that
for $d=2$, $q\geq2$ and $p<p_c(q,d)$, if $a_i\leq1/2$ for all $i\in
S$, then (\ref{noperc}) holds. Using the main result in \cite{KW},
this can be extended to $d=2$, $q\geq1$ and $p<p_c(q,d)$ with the same proof.
Combining this with Proposition \ref{almquas}, we obtain almost sure
quasilocality when $d=2$ for these parameters.
\begin{Corollary}
If $q\geq1, p<p_c(q,2)$ and $a_i\leq1/2$ for all $i\in S$, then $\mu
^{\mathbb{Z}^2}_{p,q,(a_1,\ldots,a_s)}$ is almost surely quasilocal.
\end{Corollary}

%s3 ###
\section{Useful tools}\label{usefultools}

Here, we collect the lemmas needed for the proofs of the results in
Section \ref{mainresults}.
The statements of the most important ones, Lemma \ref{condedgezd} and
Corollary \ref{barrierusezd}, are proved for finite graphs first, then
a limit is taken.
We will have an appropriate limiting procedure only for $q\geq1$ and
this is the reason why we need to restrict our attention to this case
in all our results.
Throughout this section and the next, we will use the following
notation. For a set $W\subset\mathbb{Z}^d$ and a spin configuration
$\sigma\in S^W$ on $W$, we define $C_W ^{\sigma}=\{(\xi,\eta)\in
\Omega\dvtx\xi_W=\sigma\}$. Analogously, for $E\subset\mathcal{E}^d$
and a bond configuration $\zeta\in\{0,1\}^E$ on $E$, we define $D_E
^{\zeta}=\{(\xi,\eta)\in\Omega\dvtx\eta_E=\zeta\}$.

For fixed parameters $s\in\{2,3,\ldots\}$, $p$, $a_1$, $a_2,\ldots
,a_s$ and $q\geq1$, the measure $\mathbb{P}^{\mathbb
{Z}^d}_{p,q,(a_1,\ldots,a_s)}$ can be\vspace*{1pt} obtained as a limit as follows.
Let $G_n=(\mathcal{V}_n,\mathcal{E}_n)$ be as in Section \ref{rcmeasures}.
Consider the $\operatorname{DaC}(q)$ model on $G_n$ with the given parameters as
defined in the \hyperref[intro]{Introduction}.
The corresponding sequence of measures $\mathbb
{P}^{G_n}_{p,q,(a_1,\ldots,a_s)}$ then converges to $\mathbb
{P}^{\mathbb{Z}^d}_{p,q,(a_1,\ldots,a_s)}$ as $n\to\infty$, in the
sense that probabilities of cylinder sets converge. Note that $q\geq1$
is needed to ensure the convergence of $\Phi_{p,q}^{G_n}$ to the
(unique) random-cluster measure
$\Phi_{p,q}^{\mathbb{Z}^d}$; see Section \ref{rcmeasures}.

The next two lemmas, which give the conditional edge distribution in
the $\operatorname{DaC}(q)$ model given any spin configuration, are of crucial
importance for the rest of this paper. The statements (and the proofs)
are analogs of Proposition 5.1 and Theorem 6.2 in \cite{fuzzyGibbs}.
For a graph $G=(\mathcal{V},\mathcal{E})$ (where $\mathcal{V}$ and
$\mathcal{E}$ are finite or $\mathcal{V}=\mathbb{Z}^d$, $\mathcal
{E}=\mathcal{E}^d$) and a spin configuration $\sigma\in\Omega_C^G$,
we define, for all $i\in S$,
the vertex sets $\mathcal{V}^{\sigma,i}=\{v\in\mathcal{V}\dvtx\sigma
(v)=i\}$, edge sets $\mathcal{E}^{\sigma,i}=\{e= \langle x,y
\rangle\dvtx x,y\in\mathcal{V}^{\sigma,i}\}$, $\mathcal{E}^{\sigma,
\mathrm{diff}}=\mathcal{E}\setminus\bigcup_{i=1}^s\mathcal{E}^{\sigma,i}$ and
graphs $G^{\sigma,i}=(\mathcal{V}^{\sigma,i},\mathcal{E}^{\sigma,i})$.
\begin{Lemma}\label{condedgefin}
Let $G=(\mathcal{V},\mathcal{E})$ be a finite graph. Fix parameters
$p\in[0,1]$, $q>0$, $s\in\{2,3,\ldots\},a_1,a_2,\ldots,a_s\in
(0,1)$ such that $\sum_{i=1}^{s}a_i=1$ and an arbitrary spin
configuration $\sigma\in S^{\mathcal{V}}$, and define the event $A=\{
(\xi,\eta)\in\Omega^G\dvtx \xi=\sigma\}$.
We then have have that:
\begin{itemize}
\item[(a)] for all $e\in\mathcal{E}^{\sigma,\mathrm{diff}}$, $\mathbb
{P}^G_{p,q,(a_1,\ldots,a_s)}(\{(\xi,\eta)\in\Omega^G\dvtx \eta(e)=0\}
\mid A)=1$;
\item[(b)] for all $i\in S$, independently for different values of $i$,
on the set $\{0,1\}^{\mathcal{E}^{\sigma,i}}$, the conditional
distribution of $\mathbb{P}^G_{p,q,(a_1,\ldots,a_s)}$ given $A$ is the
random-cluster measure $\Phi^{G^{\sigma,i}} _{p,qa_i}$.
\end{itemize}
\end{Lemma}
\begin{pf}
Statement (a) is immediate from the definition of the model.
Now, let $\eta\in\Omega_D^G$
be such that $\eta(e)=0$ for all $e\in\mathcal{E}^{\sigma,\mathrm{diff}}$.
Denote by $k^{\sigma,i}(\eta)$ the number of connected components in
$\eta$ that have spin $i$ in $\sigma$
and note that $k(\eta)=\sum_{i=1}^{s}k^{\sigma,i}(\eta)$. Using
this observation, (\ref{FKfinite}) and a rearrangement of the factors,
we obtain that
\begin{eqnarray*}
\mathbb{P}
((\sigma,\eta)) & = & \Phi^{G}_{p,q} ( \eta ) \prod
_{i=1}^{s}a_i^{k^{\sigma,i}(\eta)}\\
& = & \frac{(1-p)^{|\mathcal{E}^{\sigma,\mathrm{diff}}|}}{Z^G_{p,q}}\prod_{i=1}^s \biggl( (qa_i)^{k^{\sigma,i}(\eta
)}\prod_{e\in\mathcal{E}^{\sigma,i}}p^{\eta(e)}(1-p)^{1-\eta
(e)} \biggr),
\end{eqnarray*}
where we have written $\mathbb{P}$ for $\mathbb{P}^{G}_{p,q,(a_1,\ldots
,a_s)}$ and $|\cdot|$ for cardinality. It follows that
\[
\mathbb{P}^{G}_{p,q,(a_1,\ldots,a_s)}((\sigma,\eta)\mid A)=\prod
_{i=1}^{s}\Phi_{p,qa_i}^{G^{\sigma,i}}(\eta_{\mathcal{E}^{\sigma,i}})
\]
since the\vspace*{1pt} factor
$\frac{(1-p)^{|\mathcal{E}^{\sigma,\mathrm{diff}}|}\prod_{i=1}^s{Z}_{p,qa_i}^{G^{\sigma,i}}}{Z^G_{p,q}\mu
^{G}_{p,q,(a_1,\ldots,a_s)}(\sigma)}$ is constant in $\eta$ and,
thus, it must be 1 to give a probability measure. This proves statement (b).
\end{pf}
\begin{Remark}\label{condedgerem}
Let $\sigma\in S^{\mathcal{V}}$ and $A\subset\Omega^G$ be as in
Lemma \ref{condedgefin}.
The fact that random-cluster measures factorize on disconnected graphs
provides a simple way of drawing a random bond configuration $Y$ with
distribution $\mathbb{P}^G_{p,q,(a_1,\ldots,a_s)}$ given $A$.
First, set $Y(e)=0$ for all $e\in\mathcal{E}^{\sigma,\mathrm{diff}}$.
Then choose any component $C=(\mathcal{V}_C,\mathcal{E}_C)$
in the graph $(\mathcal{V},\mathcal{E}\setminus\mathcal{E}^{\sigma
,\mathrm{diff}})$.
Note that $C$ is a maximal monochromatic component in $G$ (with respect
to $\sigma$);
suppose that for all $v\in\mathcal{V}_C$, $\sigma(v)=i$.
Then, independently of everything else, draw $Y _{\mathcal{E}_C}$
according to the random-cluster measure $\Phi^{C}_{p,qa_i}$.
Repeat this procedure with a new component in $(\mathcal{V},\mathcal
{E}\setminus\mathcal{E}^{\sigma,\mathrm{diff}})$
until there are no more such components. Lemma \ref{condedgefin} and
the observation at the beginning of this paragraph ensure that we get
the correct (conditional) distribution.
\end{Remark}

By using Lemma \ref{condedgefin} and the limiting procedure for $\mathbb
{P}^{\mathbb{Z}^d}_{p,q,(a_1,\ldots,a_s)}$, one obtains analogous
statements for $\mathbb{Z}^d$ in the case $q\geq1$.
\begin{Lemma}\label{condedgezd}
Fix parameters $d$, $p$, $q\geq1$, $s$, $(a_1,a_2,\ldots,a_s)$ of the
$\operatorname{DaC}(q)$ model on $\mathbb{Z}^d$ and a spin configuration $\sigma\in
\Omega_C$.
The conditional distribution of $\mathbb{P}=\mathbb{P}^{\mathbb
{Z}^d}_{p,q,(a_1,\ldots,a_s)}$
given $C_{\mathbb{Z}^d}^{\sigma}$
then assigns value $0$ to all edges in $\mathcal{E}^{\sigma,\mathrm{diff}}$
and is a random-cluster measure for $G^{\sigma,i}$ with
parameters $p$ and $qa_i$ on $\mathcal{E}^{\sigma,i}$, independently
for each $i$. Moreover, for each edge $e\in\mathcal{E}^{\sigma,i}$
and almost every edge configuration $\zeta\in\{0,1\}^{\mathcal
{E}^d\setminus\{e\}}$,
we have that
\[
\mathbb{P}\bigl(\{(\xi,\eta)\in\Omega\dvtx\eta(e)=1\} \mid C_{\mathbb
{Z}^d}^{\sigma}\cap D_{\mathcal{E}^d\setminus\{e\}}^{\zeta}\bigr)=
\cases{
p, &\quad if $\stackrel{\zeta}{x\leftrightarrow y}$,\cr
\dfrac{p}{p+(1-p)qa_i}, &\quad otherwise.}
\]
\end{Lemma}
\begin{pf*}{Proof sketch} Unless the edge configuration $\zeta\in
\{0,1\}^{\mathcal{E}^d\setminus\{e\}}$
is special, in the sense that it contains at least two infinite FK clusters
or there exists an edge $f\in\mathcal{E}^d\setminus\{e\}$ such that
changing the state of $f$ in $\zeta$ would create at least two
infinite FK clusters,
we see,\vspace*{-1pt} after a certain stage of the limiting construction described
at the beginning of this section, whether or not
$\stackrel{\zeta}{x\leftrightarrow y}$ occurs; therefore, an equality
corresponding to the ``moreover'' part of Lemma \ref{condedgezd} can be
verified by Lemma \ref{condedgefin} for all further stages of the
limiting construction. Since the aforementioned special edge configurations
have $\Phi^{\mathbb{Z}^d}_{p,q}$-measure $0$, we are done. For the
details, see the proof of Theorem 6.2 in \cite{fuzzyGibbs}.
\end{pf*}

The next lemma, which is a more general form of Lemma 7.3 in \cite
{fuzzyGibbs}, and which can be proven in the same way, shows that,
given edge and spin configurations of a certain type [such as the ones
that we shall use in the proof of Theorem \ref{mainthm} parts 1, 2(b),
and 3(b);
see Figure \ref{configurations} before Lemma \ref{keylemma}],
the ``price of changing a spin'' depends only on the existence or
nonexistence of connections in the edge configuration. Since it appears
somewhat specialized and will not be used until Section \ref
{mainproofs}, the reader might choose to skip it for now.
\begin{Lemma}\label{price-of-change}
Fix parameters $d\geq2$, $q\geq1$, $p\in[0,1)$, $s$ and
$(a_1,a_2,\ldots,a_s)$ of the $\operatorname{DaC}(q)$ model
and let $i,j\in S$ be different spin values.
There then exist positive constants $c_1^{i,j}=c_1^{i,j}(p,q,a_i,a_j)$
and $c_2^{i,j}=c_2^{i,j}(p,q,a_i,a_j)$ such that for any $v\in\mathbb
{Z}^d$ with nearest neighbors $u_1,u_2,\ldots,u_{2d}$ and the edges
between $v$ and $u_i$ denoted by $e_i$ ($i\in\{1,2,\ldots,2d\}$), we
have, for all $\sigma\in S^{\mathbb{Z}^d\setminus\{v\}}$ and $\zeta
\in\{0,1\}^{\mathcal{E}^d\setminus\{e_1,e_2,\ldots,e_{2d}\}}$
satisfying:
\begin{enumerate}
\item$\sigma(u_1)=\sigma(u_2)=i$ and $\sigma(u_3)=\sigma
(u_4)=\cdots=\sigma(u_{2d})=j$, and
\item no two of $u_3,u_4,\ldots,u_{2d}$ are connected in $\zeta$,
\end{enumerate}
that
\[
\frac{\mathbb{P}^{\mathbb{Z}^d}_{p,q,(a_1,\ldots,a_s)}(C_{\{ v\}
}^i\mid C_{\mathbb{Z}^d\setminus\{v\}}^{\sigma}\cap D_{\mathcal
{E}^d\setminus\{e_1,\ldots,e_{2d}\}}^{\zeta})}{\mathbb{P}^{\mathbb
{Z}^d}_{p,q,(a_1,\ldots,a_s)}(C_{\{ v\}}^j\mid C_{\mathbb
{Z}^d\setminus\{v\}}^{\sigma}\cap D_{\mathcal{E}^d\setminus\{
e_1,\ldots,e_{2d}\}}^{\zeta})}=
\cases{
c_1^{i,j}, &\quad if $\stackrel{\zeta}{u_1\leftrightarrow u_2}$,\vspace*{1pt}\cr
c_2^{i,j}, &\quad otherwise.}
\]
The exact values of $c_1^{i,j}$ and $c_2^{i,j}$ are
\[
c_1^{i,j}=\frac
{p^2qa_i+2p(1-p)qa_i+(1-p)^2(qa_i)^2}{(1-p)^2(qa_i)^2}\cdot\frac
{a_i}{a_j}\cdot \biggl( \frac{(1-p)qa_j}{p+(1-p)qa_j} \biggr) ^{2d-2}
\]
and
\[
c_2^{i,j}=\frac{p^2+2p(1-p)qa_i+(1-p)^2(qa_i)^2}{(1-p)^2(qa_i)^2}\cdot
\frac{a_i}{a_j}\cdot \biggl( \frac{(1-p)qa_j}{p+(1-p)qa_j} \biggr) ^{2d-2},
\]
and this shows that
\[
\cases{
c_1^{i,j}>c_2^{i,j}& if and only if \quad$qa_i>1$,\vspace*{1pt}\cr
c_1^{i,j}=c_2^{i,j} & if and only if \quad$qa_i=1$,\vspace*{1pt}\cr
c_1^{i,j}<c_2^{i,j} & if and only if \quad$qa_i<1$.}
\]
\end{Lemma}

Lemma \ref{price-of-change} will play a role in proving parts 1, 2(b),
and 3(b) of Theorem \ref{mainthm}.
For the proof of parts 2(a) and 3(a), we will need Lemma \ref{stochlarger},
which is preceded by a few definitions and another lemma. The next definition
is motivated by Corollary~\ref{barrierusezd}.
\begin{Definition}\label{barrierdef}
We call an edge set $E=\{e_1,e_2,\ldots,e_k\}$ a \textit{barrier} if
removing $e_1,e_2,\ldots,e_k$ (but not their endvertices) separates
the graph $\mathbb{Z}^d$ into two
or more disjoint connected subgraphs.
Note that exactly one of the resulting subgraphs is infinite; we call
this the \textit{exterior} of $E$ and denote it by $\operatorname{ext}(E)$.
We denote the vertex set of $\operatorname{ext}(E)$ by $\mathcal{V}_{\mathrm{ext}(E)}$ and the
edge set of $\operatorname{ext}(E)$ by $\mathcal{E}_{\mathrm{ext}(E)}$.
We call the union of the finite subgraphs the \textit{interior} of $E$
and denote it by $\operatorname{int}(E)$. We use
$\mathcal{V}_{\mathrm{int}(E)}$ and $\mathcal{E}_{\mathrm{int}(E)}$ to denote its vertex
and edge set, respectively.
$E=\{e_1,e_2,\ldots,e_k\}$ is called a \textit{closed barrier} in a
configuration $(\xi,\eta)\in\Omega$
if $E$ is a barrier and $\eta(e_i)=0$ holds for all $i\in\{1,2,\ldots
,k\}$ and it is called a \textit{quasi-closed barrier} if,
for all edges $e= \langle x,y \rangle \in E$ such that $\eta(e)=1$, it
holds that $\xi(x)=\xi(y)\in S_{1/q}$.
\end{Definition}

For a vertex set $W \subset\mathbb{Z}^d$, we define the \textit{edge
boundary} $\Delta W$ of $W$ by
$\Delta W=\{ \langle x,y \rangle  \in\mathcal{E}^d \dvtx x\in W, y\in\mathbb
{Z}^d \setminus W\}$.
Note that the edge boundary of a union of finite spin clusters is a
closed barrier and that all closed barriers are quasi-closed.

According to Lemma \ref{condedgezd}, the states of edges in spin $i$
clusters where $a_i=1/q$ are chosen
independently of everything else, hence they should play no role in
issues of dependence.
We prove a formal statement concerning this in the following lemma and
in Corollary \ref{barrierusezd}, we show
a way to make use of this feature of the model.
For an event $A$, we denote the indicator random variable of $A$ by
$\mathbb{I}_A$.
\begin{Lemma}\label{barrieruse}
Let $G=(\mathcal{V},\mathcal{E})$ be a finite graph, let
$V_1,V_2\subset\mathcal{V}$ be a partition of $\mathcal{V}$
and, for $i\in\{1,2\}$, define edge sets $E_i=\{e\in\mathcal{E}$:
both endvertices of $e$ are in $V_i\}$
and graphs $G_i=(V_i,E_i)$. Also, define the edge set
$B=\{e\in\mathcal{E}\dvtx e$ has one endvertex in $V_1$ and one in $V_2\}$,
a subset $B^0\subset B$ and, for $i\in\{1,2\}$, $W_i$ as the set of
endvertices of edges in $B\setminus B^0$
that are in $V_i$; see Figure \ref{barrierfigure}.
%
%f1 ###
\begin{figure}

\includegraphics{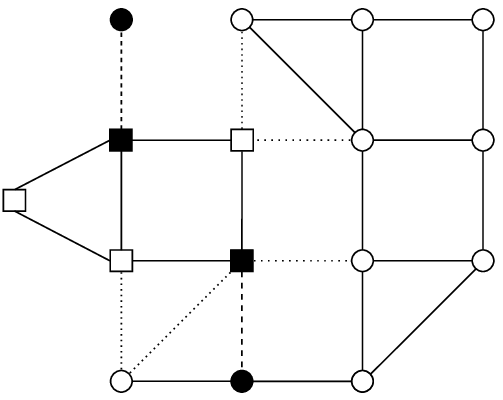}

\caption{Illustration of the situation considered in Lemma
\protect\ref{barrieruse}. The circles represent the vertices in $V_1$
and the
squares represent the vertices in $V_2$. The union of the dotted and
dashed edges makes up~$B$, with the dotted ones being in $B_0$.
Accordingly, the black circles represent the vertices in $W_1$ and the
black squares represent the vertices in $W_2$.}
\label{barrierfigure}\vspace*{-3pt}
\end{figure}
Fix parameters $p,q>0$, $s$, $(a_1,a_2,\ldots,a_s)$ of the $\operatorname{DaC}(q)$
model on $G$
and a spin configuration $\sigma\in S^{W_1\cup W_2}_{1/q}$ such that
for all $e= \langle x,y \rangle \in B\setminus B^0$,
we have that $\sigma(x)=\sigma(y)$. Considering the events
$C(B,\sigma)=\{(\xi,\eta)\in\Omega^G\dvtx \eta_{B^0}\equiv0, \xi
_{W_1\cup W_2}=\sigma\}$,
$K_1=\{(\xi,\eta)\in\Omega^{G_1}\dvtx\xi_{W_1}=\sigma_{W_1}\}$,
$K_2=\{(\xi,\eta)\in\Omega^{G_2}\dvtx\xi_{W_2}=\sigma_{W_2}\}$ and
$Z(B^0)=\{\eta\in\{0,1\}^{B^0}\dvtx\eta\equiv0\}$,
we have, for each $(\xi,\eta)\in\Omega^G$, that
\begin{eqnarray*}
&&
\mathbb{P}^G_{p,q,(a_1,a_2,\ldots,a_s)}((\xi,\eta)\mid C(B,\sigma))\\
&&\qquad = \mathbb{P}^{G_1}_{p,q,(a_1,a_2,\ldots,a_s)}((\xi_{V_1},\eta
_{E_1})\mid K_1) \times\mathbb{P}^{G_2}_{p,q,(a_1,a_2,\ldots,a_s)}((\xi
_{V_2},\eta_{E_2})\mid K_2)\nonumber\\
&&\qquad\quad{} \times\mathbb{I}_{Z(B^0)}(\eta_{B^0})\prod_{e\in B\setminus
B^0}p^{\eta(e)}(1-p)^{1-\eta(e)}.
\end{eqnarray*}
This implies, in particular, the conditional independence given
$C(B,\sigma)$ of the random configurations on $G_1$ and on $G_2$.
\end{Lemma}
\begin{pf}
Let us fix $(\xi,\eta)\in\Omega^G$. Note that
\[
\mathbb{I}_{C(B,\sigma)}(\xi,\eta)=\mathbb{I}_{K_1}(\xi_{V_1},\eta
_{E_1})\mathbb{I}_{K_2}(\xi_{V_2},\eta_{E_2})\mathbb{I}_{Z(B_0)}(\eta_{B_0}).
\]
Hence, if $(\xi,\eta)\notin C(B,\sigma)$, then we have that both
sides of the equation that we want to prove are $0$; thus, for all such
configurations, we indeed have equality of the two sides. Therefore,
let us assume that $(\xi,\eta)\in C(B,\sigma)$.
Define the event $A=\{(\kappa,\zeta)\in\Omega^G$: there is no edge
$e= \langle x,y \rangle  \in\mathcal{E}$ with $\zeta(e)=1$ and
$\kappa(x)\neq\kappa(y)\}$
and denote the analogously defined subsets of $\Omega^{G_1}$ and
$\Omega^{G_2}$ by $A_1$ and $A_2$, respectively.
Since $(\xi,\eta)\in C(B,\sigma)$, we have that
\[
\mathbb{I}_A(\xi,\eta)=\mathbb{I}_{A_1}(\xi_{V_1},\eta_{E_1})\mathbb
{I}_{A_2}(\xi_{V_2},\eta_{E_2}).
\]
Therefore, if $(\xi,\eta)\notin A$, we have $0$ on both sides of the
desired equation in Lem\-ma~\ref{barrieruse}, by the definition of the
model, so let us assume that
$(\xi,\eta)\in A$.

Now, denote by $n$ the total number of FK clusters in $\eta$. For all
$i\in S$, $j\in\{1,2\}$, denote the number of FK clusters in $\eta$
that contain
a vertex in $V_j$ with spin $i$ in $\xi$ but no vertex in $W_j$ by $n_j^i$.
For each $i\in S$, denote the number of FK clusters in $\eta$ that
contain a vertex in $W_1\cup W_2$ with spin $i$ in $\xi$ by $n_3^i$.
Throughout this proof, we shall omit the subscripts
of the joint measures in the $\operatorname{DaC}(q)$ models; for example, we write
$\mathbb{P}^G$ for the measure $\mathbb{P}^G_{p,q,(a_1,a_2,\ldots,a_s)}$.
Since $(\xi,\eta)\in C(B,\sigma)\cap A$, it immediately follows from
the definition of $\mathbb{P}^{G}$
and the definition (\ref{FKfinite}) of random-cluster measures that
%
%e5 ###
\begin{equation}\label{induncond}
\mathbb{P}^{G}((\xi,\eta)) = \frac{q^n}{Z^G _{p,q}}
\biggl(\prod_{e\in\mathcal{E}}p^{\eta(e)}(1-p)^{1-\eta(e)}\biggr)\Biggl(\prod
_{i=1}^sa_i^{n_1^i+n_2^i+n_3^i}\Biggr).
\end{equation}
Note that $\mathcal{E}=E_1\cup E_2\cup B^0\cup(B\setminus B^0)$. Since
$(\xi,\eta)\in C(B,\sigma)\cap A$, we have that
$\prod_{e\in B^0}p^{\eta(e)}(1-p)^{1-\eta(e)}=(1-p)^{|B^0|}$, where
$|\cdot|$ denotes cardinality, and
$n=\sum_{i=1}^sn_1^i+n_2^i+n_3^i$. Furthermore, it is the case that
\[
\prod_{i=1}^s(qa_i)^{n_1^i+n_2^i+n_3^i}=\prod_{i=1}^s(qa_i)^{n_1^i+n_2^i}
\]
since for all $i\notin S_{1/q}$, we have $n_3^i=0$, whereas for all
$i\in S_{1/q}$,
we have $qa_i=1$, so the factor $\prod_{i=1}^s(qa_i)^{n_3^i}$ is
indeed $1$.
Using these observations,
we can factorize the expression in (\ref{induncond}).
Indeed, denoting by $c$ the quantity $(1-p)^{|B^0|}/(Z^G_{p,q}\mathbb
{P}^G(C(B,\sigma)))$ which does
not depend on $(\xi,\eta)$, we have that
\begin{eqnarray*}
\mathbb{P}^{G}((\xi,\eta)\mid C(B,\sigma)) & = & \frac{\mathbb
{P}^{G}((\xi,\eta))}{\mathbb{P}^G(C(B,\sigma))} \\
& = & c \Biggl[ \prod_{e\in E_1}p^{\eta(e)}(1-p)^{1-\eta(e)}\prod
_{i=1}^s(qa_i)^{n_1^i} \Biggr]\\
& &{} \times \Biggl[ \prod_{e\in E_2}p^{\eta(e)}(1-p)^{1-\eta(e)}\prod
_{i=1}^s(qa_i)^{n_2^i} \Biggr] \\
& &{} \times \biggl[ \prod_{e\in B\setminus B_0}p^{\eta(e)}(1-p)^{1-\eta
(e)} \biggr].
\end{eqnarray*}

The last part of the proof, that is, showing that the expressions
between the first and second pairs of square brackets are
$c_1\mathbb{P}^{G_1}((\xi_{V_1},\eta_{E_1})\mid K_1)$ and $c_2\mathbb
{P}^{G_2}((\xi_{V_2},\eta_{E_2})\mid K_2)$, respectively,
where $c_1$ and $c_2$ are constants (i.e., they do not depend on $\xi$
or $\eta$) will be easy.
It is sufficient to show the first of these since the second one then
follows by relabeling $V_1$ and $V_2$.
Let $n_4$ denote the total number of FK clusters in $\eta_{E_1}$ and,
for each $i\in S$, let $n_5^i$ denote the number of FK clusters in
$\eta_{E_1}$
that contain a vertex in $W_1$ with spin $i$ in $\xi_{V_1}$. Since
$(\xi,\eta)\in C(B,\sigma)\cap A$, we have that $n_4=\sum
_{i=1}^sn_1^i+n_5^i$. Similarly as in the paragraph after (\ref{induncond}),
we have $n_5^i=0$ for all $i\notin S_{1/q}$ and $qa_i=1$ for all $i\in
S_{1/q}$, so
\[
\prod_{i=1}^s(qa_i)^{n_1^i+n_5^i}=\prod_{i=1}^s(qa_i)^{n_1^i}.
\]
Denoting $Z^{G_1}_{p,q}\mathbb{P}^{G_1}(K_1)$ by $c_1$, the above
observations imply that
\begin{eqnarray*}
\mathbb{P}^{G_1}((\xi_{V_1},\eta_{E_1})\mid K_1) & = & \frac{\mathbb
{P}^{G_1}((\xi_{V_1},\eta_{E_1}))}{\mathbb{P}^{G_1}(K_1)} \\
& = & \frac{q^{n_4}}{c_1} \prod_{e\in E_1}p^{\eta(e)}(1-p)^{1-\eta
(e)}\prod_{i=1}^sa_i^{n_1^i+n_5^i}\\
& = & \frac{1}{c_1} \prod_{e\in E_1}p^{\eta(e)}(1-p)^{1-\eta
(e)}\prod_{i=1}^s(qa_i)^{n_1^i}.
\end{eqnarray*}

Finally, observe that none of $c,c_1,c_2$ depend on $(\xi,\eta)$,
hence the product $cc_1c_2$ must be equal to $1$
to make $\mathbb{P}^G(\cdot\mid C(B,\sigma))$ a probability measure.
This observation completes the proof of Lemma \ref{barrieruse}.
\end{pf}

Lemma \ref{barrieruse},\vspace*{1pt} combined with the limiting procedure for
$\mathbb{P}^{\mathbb{Z}^d}_{p,q,(a_1,\ldots,a_s)}$, yields the
following result, which shows why quasi-closed barriers are useful.
\begin{Corollary}\label{barrierusezd}
Fix parameters $d$, $p$, $q\geq1$, $s$ and $(a_1,a_2,\ldots,a_s)$ of
the $\operatorname{DaC}(q)$ model on $\mathbb{Z}^d$.
Let $(X,Y)$ be a random configuration in $\Omega$ with distribution
$\mathbb{P}^{\mathbb{Z}^d}_{p,q,(a_1,\ldots,a_s)}$, $B$ a barrier\vspace*{2pt} and
$C(B)$ the event that $B$ is quasi-closed.
Then, given $C(B)$, $(X_{\mathcal{V}_{\mathrm{int}(B)}},Y_{\mathcal
{E}_{\mathrm{int}(B)}})$ and $(X_{\mathcal{V}_{\mathrm{ext}(B)}},Y_{\mathcal
{E}_{\mathrm{ext}(B)}})$ are conditionally independent.
In particular, for
a set $H\subset\mathbb{Z}^d$ and a spin configuration $\sigma\in S^H$,
we have that the conditional distribution of $(X_{\mathcal
{V}_{\mathrm{int}(B)}},Y_{\mathcal{E}_{\mathrm{int}(B)}})$ given $C(B)$ and
$\{(\xi,\eta)\in\Omega\dvtx \xi_H=\sigma_H\}$
is $\mathbb{P}^{\mathrm{int}(B)}_{p,q,(a_1,\ldots,a_s)}$ conditioned on $\{(\xi
,\eta)\in\Omega^{\mathrm{int}(B)}\dvtx \xi_{H\cap\mathcal{V}_{\mathrm{int}(B)}}=\sigma
_{H\cap\mathcal{V}_{\mathrm{int}(B)}}\}$.
\end{Corollary}
\begin{Remark}\label{PottsMarkovian}
As a first application of Corollary \ref{barrierusezd}, we give a
proof of the Markovianness of the measure $\mu^{\mathbb
{Z}^d}_{p,q,(1/q,1/q,\ldots,1/q)}$ that does not use the connection
between the $\operatorname{DaC}(q)$ and Potts models.
Fix a finite subset $W$ of $\mathbb{Z}^d$.
For any spin configuration $\sigma\in S^{\mathbb{Z}^d\setminus W}$
outside $W$, we have that
the edge set $B=\{e\in\mathcal{E}^d\dvtx e$ has one end-vertex in $\partial
W$ and one in $\partial_2W\setminus\partial W\}$ is a quasi-closed
barrier since,
for each edge $e= \langle x,y \rangle  \in B$, it is either the case that:
$\sigma(x)\neq\sigma(y)$ and, therefore, $e$ is closed;
or, $\sigma(x)=\sigma(y)\in S_{1/q}$ since $S=S_{1/q}$. Therefore, we
have, by Corollary~\ref{barrierusezd}, that for any spin configurations
$\sigma,\sigma^{\prime}\in S^{\mathbb{Z}^d\setminus W}$ with $\sigma
^{\prime}_{\partial W}=\sigma_{\partial W}$, the conditional distributions
$\mathbb{P}^{\mathbb{Z}^d}_{p,q,(1/q,1/q,\ldots,1/q)}$ given
$K_{\mathbb{Z}^d\setminus W}^{\sigma}$ and
$\mathbb{P}^{\mathbb{Z}^d}_{p,q,(1/q,1/q,\ldots,1/q)}$ given
$K_{\mathbb{Z}^d\setminus W}^{\sigma^{\prime}}$ are the same in
$S^{\mathcal{V}_{\mathrm{int}(B)}}\times\{0,1\}^{\mathcal{E}_{\mathrm{int}(B)}}$. The
statement follows.
\end{Remark}

Corollary \ref{barrierusezd} enables us to describe some of the
intuition behind our main results.
Roughly speaking, quasilocality
means that, conditioning on a spin configuration $\sigma\in S^{\mathbb
{Z}^d\setminus H}$ outside a set $H$,
the spin distribution in $H$ does not depend on spins very far away
from $H$.
Corollary \ref{barrierusezd} shows that this is the case if $H$ is
surrounded by a quasi-closed barrier.
In particular,
the presence of such a quasi-closed barrier is automatic if
$S=S_{1/q}$, as was noted in Remark \ref{PottsMarkovian}. It is also
easy to see that if there is no percolation in $\sigma$ in any spin
(i.e., there exists no infinite connected component of equal spins),
then there exists a closed barrier surrounding $H$, namely the edge
boundary of the (finite) union of $H$ and the spin clusters in $\sigma
$ that contain at least one vertex in $\partial H$. This reasoning will
be used in the proof of Proposition \ref{almquas}.

However, if there exists an infinite spin $i$ cluster $C$ in $\sigma$
with $i\in S\setminus S_{1/q}$ that contains a vertex in $\partial H$,
then it cannot be decided whether or not a quasi-closed barrier
surrounding $H$ exists just by looking at the spin configuration; one
also needs to check the edge configuration in $C$. Clearly, if we see
an infinite open edge component in $C$ that contains a vertex in
$\partial H$, then there is no quasi-closed barrier that surrounds
$H\cup\partial H$. Since, by Lemma \ref{condedgezd}, the conditional
edge distribution in the spin $i$ cluster $C$ is a random-cluster
measure with parameters $p$ and $qa_i$, the question is whether such
measures percolate.

Now, recall the definition of $\ell\in S$ which was given immediately
before Theorem \ref{mainthm} and consider the case when $a_{\ell}>
1/q$. Note that the condition $p<p_c(qa_{\ell},d)$, which appears in
part 2(a) of Theorem \ref{mainthm}, ensures that there is no infinite
edge cluster in any spin $j$ cluster where $j\in S\setminus S_{1/q}$,
by the definition of $\ell$ and using the well-known fact (see, e.g.,
\cite{grimmett2}) that if $q_1\geq q_2$, then $p_c(q_1,d)\geq
p_c(q_2,d)$. In Section \ref{mainproofs}, we will show that for all
such $p$, there exists a quasi-closed barrier surrounding $H$ given any
spin configuration $\sigma\in S^{\mathbb{Z}^d\setminus H}$ with
arbitrarily high probability and, hence, quasilocality holds.

This argument suggests that the best candidate for a spin configuration
in which spins arbitrarily far away from $H$ still have a significant
influence on the spin distribution in $H$ (thereby implying
nonquasilocality) are those with an infinite spin $\ell$ cluster, and
that quasilocality might fail for all $p>p_c(qa_{\ell},d)$. We have
not managed to prove this, but we\vspace*{1pt} get very close by proving
nonquasilocality for all $p>p_c^{\mathcal{H}}(qa_{\ell},d)$ in
Section \ref{mainproofs}. Indeed, this is equivalent to the full
statement under the widely accepted conjecture that for random-cluster measures,
the critical value and the half-space critical value coincide.

The case where $a_{\ell}<1/q$ is more problematic, mainly because the
random-cluster measures with parameters $p$ and $qa_{\ell}$ then do
not have the property of positive association; in fact, it is
easy to see that for such random-cluster measures, the conditional
probability given in Definition \ref{FKdef} is nonincreasing in~$\zeta$.
Since, for $a_{\ell}>1/q$, positive association plays a role both in
proving quasilocality
for small $p$ and in proving nonquasilocality
for large $p$, without it, we must resort to comparing the conditional
edge configuration given a spin configuration
to Bernoulli bond percolation. This, however, yields
worse upper (resp., lower) bounds for the quasilocality
(resp., nonquasilocality) regime, leaving a gap between the bounds.

The following definition will help us to localise quasi-closed barriers.
\begin{Definition}\label{Phatdef}
Let $(X,Y)$ be an $S^{\mathbb{Z}^d}\times\{0,1\}^{\mathcal
{E}^d}$-valued random pair with distribution $\mathbb{P}^{\mathbb
{Z}^d}_{p,q,(a_1,\ldots,a_s)}$.
Given that $(X,Y)=(\xi,\eta)$ for some $(\xi,\eta)\in\Omega$, let
$\hat Y \in\{0,1\}^{\mathcal{E}^d}$ be defined by
setting, for each $e= \langle x,y \rangle  \in\mathcal{E}^d$,
\[
\hat Y (e)=
\cases{
0, &\quad if $\xi(x)=\xi(y)\in S_{1/q}$,\cr
\eta(e), &\quad otherwise.}
\]
We write $\hat{\mathbb{P}}^{\mathbb{Z}^d}_{p,q,(a_1,\ldots,a_s)}$ for
the induced joint distribution of $(X,Y,\hat Y)$
on $\hat\Omega=S^{\mathbb{Z}^d}\times\{0,1\}^{\mathcal{E}^d}\times\{
0,1\}^{\mathcal{E}^d}$.
\end{Definition}

The next lemma, which is a generalization of Lemma 9.5 in \cite
{fuzzyGibbs}, compares the conditional distribution of $\hat Y$
given a spin configuration and $\hat Y$ outside a finite edge set $F$
to a random-cluster measure for $\mathbb{Z}^d$ with parameters $p$ and
$qa_{\ell}$ in the case where $a_{\ell}>1/q$,
and the conditional distribution of $Y$ given a spin configuration and
$\hat Y$ outside a finite edge set $F$
to Bernoulli bond percolation with parameter $\frac{p}{p+(1-p)qa_{\ell
}}$ in the case where $a_{\ell}<1/q$.
\begin{Lemma}\label{stochlarger}
Suppose that $q\geq1$ and $S\neq S_{1/q}$.
For any spin configuration $\sigma\in S^{\mathbb{Z}^d}$, edge set
$F\subset\subset\mathcal{E}^d$ and edge configurations $\zeta,\zeta
^{\prime}\in\{0,1\}^{\mathcal{E}^d\setminus F}$ such that $\zeta
^{\prime}\geq\zeta$,
defining $\hat A=\{(\xi,\eta,\hat\eta)\in\hat\Omega\dvtx\xi=\sigma
,\hat\eta_{\mathcal{E}^d\setminus F}=\zeta\}$ and $A^{\prime}=\{
\eta\in\Omega_D\dvtx\eta_{\mathcal{E}^d\setminus F}=\zeta^{\prime}\}$,
we have the following:
\begin{enumerate}
\item if $a_{\ell}>1/q$ and $\phi$ is a random-cluster measure for
$\mathbb{Z}^d$ with parameters $p$ and $qa_{\ell}$, then
the conditional distribution of $\phi$ given $A^{\prime}$ is
stochastically larger than
the marginal on $\hat Y$
of $\hat{\mathbb{P}}^{\mathbb{Z}^d}_{p,q,(a_1,\ldots,a_s)}$ given
$\hat A$;\vspace*{2pt}
\item if $a_{\ell}<1/q$, then
the conditional distribution of the product measure\break $\Phi^{\mathbb
{Z}^d,1}_{{p}/({p+(1-p)qa_{\ell}}),1}$ given $A^{\prime}$
is stochastically\vspace*{1pt} larger than
the marginal on $Y$
of $\hat{\mathbb{P}}^{\mathbb{Z}^d}_{p,q,(a_1,\ldots,a_s)}$ given
$\hat A$.
\end{enumerate}
\end{Lemma}
\begin{pf}
We first prove part 1. By Holley's theorem on stochastic domination
(see \cite{GHM}, Theorem 4.8), it is sufficient to prove that for all
$a\in\{0,1\}$, $e= \langle x,y \rangle  \in F$,
$\zeta_{g},\zeta_{s}\in\{0,1\}^{F\setminus\{e\} }$ such that
$\zeta_{g}\geq\zeta_{s}$, if we define $B_g=\{\eta\in
\Omega_D\dvtx\eta_{F\setminus\{e\}}=\zeta_{g}\}$ and $\hat B_s=\{ (\xi
,\eta,\hat\eta)\in\hat\Omega\dvtx\hat\eta_{F\setminus\{e\}}=\zeta
_s\}$, then we have that
%
%e6 ###
\begin{equation}\label{phiprob}
\phi\bigl(\{ \eta\in\Omega_D\dvtx\eta(e)\geq a\} \mid A^{\prime}\cap B_g \bigr)
\end{equation}
is greater than or equal to
%
%e7 ###
\begin{equation}\label{Pprob}
\hat{\mathbb{P}}^{\mathbb{Z}^d}_{p,q,(a_1,\ldots,a_s)}\bigl(\{ (\xi,\eta
,\hat\eta)\in\hat\Omega\dvtx\hat\eta(e)\geq a\} \mid\hat A\cap\hat B_s\bigr).
\end{equation}
This is obvious for $a=0$. For $a=1$, using the notation $(\eta_1;\eta
_2)$ for an edge configuration which agrees with $\eta_1$ on $\mathcal
{E}^d\setminus F$ and with $\eta_2$ on $F\setminus\{e\}$, we have, by
Definition \ref{FKdef} of random-cluster measures, that (\ref{phiprob}) equals
\[
\cases{
p, &\quad if $x\stackrel{(\zeta^{\prime};\zeta
_g)}{\leftrightarrow}y$,\cr
\dfrac{p}{p+(1-p)qa_{\ell}}, &\quad if $x\stackrel{(\zeta^{\prime
};\zeta_g)}{\nleftrightarrow}y$.}
\]
For (\ref{Pprob}), we first need to check what the spins of the
endvertices $x,y$ of $e$ are in $\sigma$. Indeed, if $\sigma(x)\neq
\sigma(y)$
or $\sigma(x)=\sigma(y)\in S_{1/q}$, then $\mbox{(\ref{Pprob})}=0$ by
Definition \ref{Phatdef}.
Let us assume that $\sigma(x)=\sigma(y)=j\notin S_{1/q}$ and denote
the maximal monochromatic component (with respect to $\sigma$) in the
graph $\mathbb{Z}^d$ which
contains $x$ by $G_x$. By Lem\-ma~\ref{condedgezd}, the conditional
distribution of $Y$ given $\sigma$ is a random-cluster measure on
$G_x$ with parameters $p$ and $qa_j$. Moreover, since $j\notin
S_{1/q}$, we have that $\hat Y$ and $Y$ agree on $G_x$. Keeping these
observations in mind,
it follows that (\ref{Pprob}) equals
\[
\cases{
0, &\quad if $\sigma(x)\neq\sigma(y)$ or $\sigma
(x)=\sigma(y)\in S_{1/q}$,\cr
p, &\quad if $\sigma(x)=\sigma(y)\notin S_{1/q}$ and
$x\stackrel{(\zeta;\zeta_s)}{\leftrightarrow}y$,\cr
\dfrac{p}{p+(1-p)qa_j}, &\quad if $\sigma(x)=\sigma(y)=j\notin
S_{1/q}$ and $x\stackrel{(\zeta;\zeta_s)}{\nleftrightarrow}y$.}
\]
Since, due to the assumption $a_{\ell}\geq1/q$, we have that
\[
p\geq\frac{p}{p+(1-p)qa_j}
\]
for all $j\in S$ and, by the definition of $a_{\ell}$, we have that
\[
\frac{p}{p+(1-p)qa_{\ell}}\geq\frac{p}{p+(1-p)qa_j}
\]
for all $j\in S\setminus S_{1/q}$, we obtain the desired result by
noting that $x\stackrel{(\zeta^{\prime};\zeta_g)}{\nleftrightarrow
}y$ implies that $x\stackrel{(\zeta;\zeta_s)}{\nleftrightarrow}y$.

Part 2 can also be proven by a direct application of Holley's theorem,
noting that,
due to the definition of $\ell$ and the assumption $qa_{\ell}<1$, we have
\[
\hspace*{63pt}
\frac{p}{p+(1-p)qa_{\ell}}\geq\max\biggl\{p,\max_{i\in S}\frac
{p}{p+(1-p)qa_{i}}\biggr\}.
\hspace*{63pt}\qed
\]
\noqed\end{pf}

Although the set $F$ had to be finite in Lemma \ref{stochlarger} so
that we could use Holley's theorem in the proof, it will not be
difficult to deduce an analogous statement corresponding to $F=\mathcal
{E}^d$; see below.
\begin{Corollary}\label{stochlargercor}
Suppose that $q\geq1$ and $S\neq S_{1/q}$, and
let $\sigma\in S^{\mathbb{Z}^d}$ be an arbitrary spin configuration.
Defining $\hat A=\{(\xi,\eta,\hat\eta)\in\hat\Omega\dvtx\xi=\sigma
\}$, we have the following:
\begin{enumerate}
\item if $a_{\ell}>1/q$, then the wired random-cluster measure $\Phi
^{\mathbb{Z}^d,1}_{p,qa_{\ell}}$ is stochastically larger than
the marginal on $\hat Y$
of $\hat{\mathbb{P}}^{\mathbb{Z}^d}_{p,q,(a_1,\ldots,a_s)}$ given
$\hat A$;
\item if $a_{\ell}<1/q$, then the product measure $\Phi^{\mathbb
{Z}^d,1}_{{p}/({p+(1-p)qa_{\ell}}),1}$
is stochastically larger than the marginal on $Y$ of
$\hat{\mathbb{P}}^{\mathbb{Z}^d}_{p,q,(a_1,\ldots,a_s)}$ given $\hat A$.
\end{enumerate}
\end{Corollary}
\begin{pf}
We only give the proof of part 1 since part 2 can be proven analogously.
Assume that $a_{\ell}>1/q$ and let $\phi$ be a random-cluster measure
for $\mathbb{Z}^d$ with parameters $p$ and $qa_{\ell}$. For $n\in\{
1,2,\ldots\}$, let $\mathcal{E}_n$ and $W_n$ be as in Section \ref
{rcmeasures} and define $\phi_n$ as $\phi$ conditioned on $W_n\cap\{
\eta\in\Omega_D\dvtx\eta_{\mathcal{E}^d\setminus\mathcal{E}_n}\equiv
1\}$. It follows from Lemma \ref{stochlarger} that for each $n$, $\phi
_n$ is stochastically larger
than the marginal on $\hat Y$
of $\hat{\mathbb{P}}^{\mathbb{Z}^d}_{p,q,(a_1,\ldots,a_s)}$ given
$\hat A$. On the other hand, $\phi_n$ coincides on $\mathcal{E}_n$
with $\phi^{G_n,1}_{p,qa_{\ell}}$ (which is defined in Section \ref
{rcmeasures}), so it converges to $\Phi_{p,qa_{\ell}}^{\mathbb
{Z}^d,1}$ as $n\to\infty$. Since stochastic domination is preserved
under weak limits, this observation completes the proof.
%\rightqed
\end{pf}

Note that if $S_{1/q}=\varnothing$, then $\hat Y$ can be replaced by $Y$
in part 1 of Lemma \ref{stochlarger} and Corollary \ref
{stochlargercor}. Also, since $\hat Y\leq Y$ by definition,
we could write $\hat Y$ instead of $Y$ in part 2.

%s4 ###
\section{Proofs of the main results}\label{mainproofs}

After all of the preparation in Section \ref{usefultools}, we are now
ready to prove our main results.
The proof of Proposition \ref{untheorem} is not difficult. In fact, one
can use the same idea as is used in the proof of Lemma 5.6 in \cite
{DaC}, namely that any vertex can be isolated (i.e., incident to closed
edges only) in the edge configuration (given any spin configuration)
with probability bounded away from 0, in which case it can be assigned
any spin in $S$, independently of everything else. A formal proof
proceeds as follows.
\begin{pf*}{Proof of Proposition \ref{untheorem}}
Fix $v\in\mathbb{Z}^d$, $m\in S$ and $\sigma\in S^{\mathbb
{Z}^d\setminus\{v\}}$, and recall the definition for $W\subset\mathbb
{Z}^d$ of the event $K_W^{\sigma}\subset\Omega_C$ at the beginning
of Section \ref{mainresults} and of the analogous event $C_W^{\sigma
}\subset\Omega$
at the beginning of Section \ref{usefultools}.
Denote by $E_v$ the event that all $2d$ edges incident to $v$ are
closed. We have that
%
%e8 ###
\begin{eqnarray}\label{probxvm}
\mu^{\mathbb{Z}^d}_{p,q,(a_1,\ldots,a_s)}\bigl(K_{\{v\}}^m\mid K_{\mathbb
{Z}^d\setminus\{v\}}^{\sigma} \bigr) & \geq& \mathbb{P}^{\mathbb
{Z}^d}_{p,q,(a_1,\ldots,a_s)}\bigl(C_{\{v\}}^m\mid C_{\mathbb{Z}^d\setminus
\{v\}}^{\sigma}\cap E_v\bigr)\nonumber\\[-8pt]\\[-8pt]
& &{} \times \mathbb{P}^{\mathbb{Z}^d}_{p,q,(a_1,\ldots
,a_s)}\bigl(E_v\mid
C_{\mathbb{Z}^d\setminus\{v\}}^\sigma\bigr).\nonumber
\end{eqnarray}
Obviously (or as a special case of Corollary \ref{barrierusezd}), we
have that the first term on the right-hand side of (\ref{probxvm}) is
$a_m$ since, given $E_v$, $v$ is assigned a spin independently of
everything else.

On the other hand,
\begin{eqnarray*}
\mathbb{P}^{\mathbb{Z}^d}_{p,q,(a_1,\ldots,a_s)}\bigl(E_v\mid C_{\mathbb
{Z}^d\setminus\{v\}}^{\sigma}\bigr) & = & \sum_{b\in S}\mathbb{P}^{\mathbb
{Z}^d}_{p,q,(a_1,\ldots,a_s)}\bigl(E_v\mid C_{\mathbb{Z}^d\setminus\{v\}
}^{\sigma}\cap C_{\{ v \}}^b\bigr)\\
& &\hspace*{12.4pt}{} \times \mathbb{P}^{\mathbb{Z}^d}_{p,q,(a_1,\ldots
,a_s)}\bigl(C_{\{v\}}^b\mid C_{\mathbb{Z}^d\setminus\{v\}}^{\sigma}\bigr).
\end{eqnarray*}
Now, whatever value $b\in S$ takes, the full spin configuration is
given in the first factor on the right-hand side, so we can apply Lemma
\ref{condedgezd}.
Under any random-cluster measure with parameters $p$ and $\tilde q>0$,
the probability of $E_v$ is bounded away from 0: a lower bound for
$\tilde q\geq1$ is $(1-p)^{2d}$, while for $\tilde q<1$ it is $(1-\frac
{p}{p+(1-p)\tilde q})^{2d}$. Since the parameter $\tilde q$ here equals
$qa_b$ for some $b$, we get the lower bound
\[
\min\biggl\{(1-p)^{2d},\biggl(1-\frac{p}{p+(1-p)q\min_{i\in S}a_i}\biggr)^{2d}\biggr\}
\]
for the first factor, which is uniform in $b$.
As the $\mathbb{P}^{\mathbb{Z}^d}_{p,q,(a_1,\ldots,a_s)}(C_{\{v\}
}^b\mid C_{\mathbb{Z}^d\setminus\{v\}}^{\sigma})$ for $b\in S$ sum to
$1$, we get the same bound for $\mathbb{P}^{\mathbb
{Z}^d}_{p,q,(a_1,\ldots,a_s)}(E_v\mid C_{\mathbb{Z}^d\setminus\{v\}
}^{\sigma})$.

Combining this with (\ref{probxvm}) and the remark thereafter, we have that
\[
\varepsilon=  \Bigl( \min_{i\in S}a_i  \Bigr)  \biggl( \min\biggl\{1-p,1-\frac
{p}{p+(1-p)q\min_{i\in S}a_i}\biggr\} \biggr) ^{2d}
\]
is a lower bound for $\mu^{\mathbb{Z}^d}_{p,q,(a_1,\ldots,a_s)}(K_{\{
v\}}^m\mid K_{\mathbb{Z}^d\setminus\{v\}}^{\sigma} )$.
Since $\varepsilon$ does not depend on $v$, $m$ or $\sigma$ and is
positive for any values of $p\in[0,1)$, $q\geq1$ and $a_1,\ldots
,a_s\in(0,1)$, we conclude that $\mu^{\mathbb{Z}^d}_{p,q,(a_1,\ldots
,a_s)}$ is uniformly nonnull for such parameters.
\end{pf*}

The proof of Theorem \ref{mainthm} consists of many parts.
For the proof of parts 1, 2(b) and 3(b),
we use a counterexample that is very similar to the one given in \cite
{DaC,fuzzyGibbs} (see also \cite{vEMSS,Kulske1,Kulske2}), defined below.
Following the definition, we give Lemma \ref{keylemma}, after which
it will not be difficult to prove parts 1, 2(b) and 3(b).
Finally, we prove parts 2(a) and 3(a).
From this point on, we assume that $S\neq S_{1/q}$. First, recall
the definitions of $\Lambda
_n$ and $\ell$ from Sections \ref{rcmeasures} and \ref{mainresults},
respectively.
Fix an arbitrary spin $m\in S$ such that $m\neq\ell$ and define an
auxiliary spin configuration $\sigma^{*}\in S^{\mathbb{Z}^d}$ by
setting, for each $x=(x_1,x_2,\ldots,x_d)\in\mathbb{Z}^d$,
\[
\sigma^{*}(x)=
\cases{
m, &\quad if $x_1=0,|x_2|+|x_3|+\cdots+|x_d|=1$\cr
&\quad or $x_1=-1,|x_2|+|x_3|+\cdots+|x_d|>1$,\cr
\ell, &\quad otherwise,}
\]
and for $k\in\{1,2,\ldots\}$, spin configurations $\sigma^{k,\ell
},\sigma^{k,m}\in S^{\mathbb{Z}^d\setminus\{\mathbf{0}\}}$ (see Figure
\ref{configurations}) by
\[
\sigma^{k,\ell}(x)=
\cases{
\ell &\quad for $x\in\mathbb{Z}^d\setminus\Lambda_k$,\cr
\sigma^{*}(x) &\quad otherwise,}
\]
and
\[
\sigma^{k,m} (x)=
\cases{
m &\quad for $x\in\Lambda_{k+1}\setminus\Lambda_k$,\cr
\sigma^{k,\ell}(x) &\quad otherwise.}
\]

Denote the two nearest neighbors of $\mathbf{0}$ in $\mathbb{Z}^d$ with
$\sigma^{*}$-spin $\ell$ by $u_1=(1,0,\break0,\ldots,0)$ and
$u_2=(-1,0,0,\ldots,0)$,
the other nearest neighbors by $u_3,u_4,\break\ldots,u_{2d}$ and
%
%
%f2 ###
\begin{figure}%[b]\vspace*{-3pt}
\[
\begin{array}{ccccccccc@{\qquad}cccccccccccccc}
\ell& \ell& \ell& \ell& \ell& \ell& \ell& \ell& \ell& & & & &
m & m & m & m & m & m & m & m & m \\
\ell& \ell& \ell& m & \ell& \ell& \ell& \ell& \ell& & & & & m &
\ell& \ell& m & \ell& \ell& \ell& \ell& m \\
\ell& \ell& \ell& m & \ell& \ell& \ell& \ell& \ell& & & & & m &
\ell& \ell& m & \ell& \ell& \ell& \ell& m \\
\ell& \ell& \ell& \ell& m & \ell& \ell& \ell& \ell& & & & & m &
\ell& \ell& \ell& m & \ell& \ell& \ell& m \\
\ell& \ell& \ell& \ell& & \ell& \ell& \ell& \ell& & & & & m &
\ell& \ell& \ell& & \ell& \ell& \ell& m \\
\ell& \ell& \ell& \ell& m & \ell& \ell& \ell& \ell& & & & & m &
\ell& \ell& \ell& m & \ell& \ell& \ell& m \\
\ell& \ell& \ell& m & \ell& \ell& \ell& \ell& \ell& & & & & m &
\ell& \ell& m & \ell& \ell& \ell& \ell& m \\
\ell& \ell& \ell& m & \ell& \ell& \ell& \ell& \ell& & & & & m &
\ell& \ell& m & \ell& \ell& \ell& \ell& m \\
\ell& \ell& \ell& \ell& \ell& \ell& \ell& \ell& \ell& & & & &
m & m & m & m & m & m & m & m & m
\end{array}
\]
\caption{Restriction of $\sigma^{3,\ell}$ (to the left) and $\sigma
^{3,m}$ (to the right) to $\Lambda_4\setminus\{\mathbf{0}\}$ in two dimensions.
For all $x\in\mathbb{Z}^2\setminus\Lambda_4$, $\sigma^{3,\ell
}(x)=\sigma^{3,m}(x)=\ell$.}
\label{configurations}
\end{figure}
for $i\in\{1,2,\ldots,2d\}$, the edges between $\mathbf{0}$ and $u_i$ by $e_i$.
Most of the work needed for the proof of parts 1, 2(b) and 3(b) of
Theorem \ref{mainthm} is contained in the following lemma.
\begin{Lemma}\label{keylemma}
Fix parameters $d\geq2$, $q\geq1$, $p$ and $a_1,a_2,\ldots,a_s\in
(0,1)$ of the $\operatorname{DaC}(q)$ model on $\mathbb{Z}^d$
in such a way that $S\neq S_{1/q}$.
Considering the events
$A=\{(\xi,\eta)\in\Omega$: there exists an open path in $\eta
_{\mathcal{E}^d\setminus\{e_1,e_2,\ldots,e_{2d}\}}$ between $u_1$ and
$u_2\}$ and
$O^{\ell,m}=\{(\xi,\eta)\in\Omega\dvtx\xi(\mathbf{0})\in\{\ell,m\}\}$,
we have the following:
\begin{enumerate}
\item if, for a fixed $k\in\{1,2,\ldots\}$, we have that
\[
\mathbb{P}^{\mathbb{Z}^d}_{p,q,(a_1,\ldots,a_s)}\bigl(A \mid O^{\ell
,m}\cap C_{\mathbb{Z}^d\setminus\{\mathbf{0}\}}^{\sigma^{k,\ell}}\bigr)>0,
\]
then $\mu^{\mathbb{Z}^d}_{p,q,(a_1,\ldots,a_s)}$ is not
$k$-Markovian;\vspace*{1pt}
\item if there exists $\gamma>0$ such that, for all $k\in\{1,2,\ldots
\}$,
\[
\mathbb{P}^{\mathbb{Z}^d}_{p,q,(a_1,\ldots,a_s)}
\bigl(A \mid O^{\ell,m}\cap C_{\mathbb{Z}^d\setminus\{\mathbf{0}\}}^{\sigma
^{k,\ell}}\bigr)>\gamma,
\]
then $\mu^{\mathbb{Z}^d}_{p,q,(a_1,\ldots,a_s)}$ is not quasilocal.
\end{enumerate}
\end{Lemma}
\begin{pf}
In order to simplify notation, in this proof, we denote $\mathbb
{P}^{\mathbb{Z}^d}_{p,q,(a_1,\ldots,a_s)}$ by $\mathbb{P}$,
$C_{\mathbb{Z}^d\setminus\{\mathbf{0}\}}^{\sigma^{k,\ell}}$ by $L=L^k$
and $C_{\mathbb{Z}^d\setminus\{\mathbf{0}\}}^{\sigma^{k,m}}$ by $M=M^k$.
The first step in the proof is to derive,
for all $k\in\{1,2,\ldots\}$, inequality (\ref{trans}).
Consider the expression
%
%e9 ###
\begin{equation}\label{expression}
\bigl| \mathbb{P}\bigl(C_{\{ \mathbf{0}\}} ^{\ell}\mid L\bigr)-\mathbb{P}\bigl(C_{\{
\mathbf{0}\}} ^{\ell} \mid M\bigr) \bigr|.
\end{equation}
Note that we have
\[
\mathbb{P}\bigl(C_{\{ \mathbf{0}\}} ^{\ell}\mid L\bigr)=\mathbb{P}\bigl(C_{\{ \mathbf{0}\}}
^{\ell}\mid O^{\ell,m}\cap L\bigr)\mathbb{P}(O^{\ell,m}\mid L)
\]
and similarly for $M$.
Using this,
we obtain, via basic algebra
[i.e., first subtracting, then adding a dummy term
$\mathbb{P}(C_{\{ \mathbf{0}\}} ^{\ell}\mid O^{\ell,m}\cap M)\mathbb
{P}(O^{\ell,m}\mid L)$
in (\ref{expression}) between the absolute values
and finally using the fact that $|a-b|\geq|a|-|b|$], that
(\ref{expression}) is greater than or equal to
\begin{eqnarray*}
&& \bigl| \mathbb{P}\bigl(C_{\{ \mathbf{0}\}} ^{\ell} \mid O^{\ell,m}\cap
L\bigr)-\mathbb{P}\bigl(C_{\{ \mathbf{0}\}} ^{\ell}\mid O^{\ell,m}\cap M\bigr)
\bigr|\mathbb{P}(O^{\ell,m}\mid L)\\
&&\qquad{} - | \mathbb{P}(O^{\ell,m}\mid L)-\mathbb{P}(O^{\ell,m}\mid
M) | \mathbb{P}\bigl(C_{\{ \mathbf{0}\}} ^{\ell}\mid O^{\ell,m}\cap
M\bigr).
\end{eqnarray*}
Since $\mathbb{P}(O^{\ell,m}\mid L)=\mathbb{P}(C_{\{ \mathbf{0}\}} ^{\ell
}\mid L)+\mathbb{P}(C_{\{ \mathbf{0}\}} ^{m} \mid L)$,
we have, by uniform nonnullness (i.e., Proposition \ref{untheorem}),
that there exists $\delta>0$ such that, uniformly in $k$,
$ | \mathbb{P}(O^{\ell,m}\mid L) | \geq\delta$.
Using this observation,
noting that $\mathbb{P}(C_{\{\mathbf{0}\}} ^{\ell}\mid O^{\ell,m}\cap
M)\leq1$,
that $\mathbb{P}(O^{\ell,m}\mid L)=\mathbb{P}(C_{\{\mathbf{0}\}} ^{\ell
}\mid L)+\mathbb{P}(C_{\{\mathbf{0}\}} ^{m}\mid L)$
(and similarly for $M$) and applying the triangle inequality yields
that (\ref{expression}) is greater than or equal to
\begin{eqnarray*}
&& \bigl| \mathbb{P}\bigl(C_{\{\mathbf{0}\}} ^{\ell} \mid O^{\ell,m}\cap
L\bigr)-\mathbb{P}\bigl(C_{\{\mathbf{0}\}} ^{\ell}\mid O^{\ell,m}\cap M\bigr)
\bigr|\delta\\
&&\qquad{} - \bigl| \mathbb{P}\bigl(C_{\{\mathbf{0}\}} ^{\ell}\mid L\bigr)-\mathbb{P}\bigl(C_{\{
\mathbf{0}\}} ^{\ell}\mid M\bigr) \bigr|\\
&&\qquad{}- \bigl| \mathbb{P}\bigl(C_{\{ \mathbf{0}\}}
^{m}\mid L\bigr) -\mathbb{P}\bigl(C_{\{ \mathbf{0}\}} ^{m}\mid M\bigr) \bigr| .
\end{eqnarray*}
After a rearrangement of the terms, this gives that
%
%e10 ###
\begin{eqnarray} \label{trans}
&& 2 \bigl| \mathbb{P}\bigl(C_{\{ \mathbf{0}\}} ^{\ell}\mid L\bigr)-\mathbb
{P}\bigl(C_{\{ \mathbf{0}\}} ^{\ell}\mid M\bigr) \bigr|+ \bigl| \mathbb{P}\bigl(C_{\{
\mathbf{0}\}} ^{m}\mid L\bigr)-\mathbb{P}\bigl(C_{\{ \mathbf{0}\}} ^{m}\mid M\bigr) \bigr|
\nonumber\\[-8pt]\\[-8pt]
&&\qquad \geq \delta \bigl| \mathbb{P}\bigl(C_{\{ \mathbf{0}\}} ^{\ell}\mid O^{\ell
,m}\cap L\bigr)-\mathbb{P}\bigl(C_{\{ \mathbf{0}\}} ^{\ell}\mid O^{\ell,m}\cap
M\bigr) \bigr|.\nonumber
\end{eqnarray}

From now on, we will be working on bounding the right-hand side of (\ref
{trans}) from below.
Elementary calculations and an application of Lemma \ref
{price-of-change} with $i=\ell$, $j=m$ and $v=\mathbf{0}$
show that
\begin{eqnarray*}
&&\mathbb{P}\bigl(C_{\{ \mathbf{0}\}} ^{\ell}\mid O^{\ell,m}\cap L\cap A\bigr) \\
&&\qquad =
\frac{\mathbb{P}(C_{\{ \mathbf{0}\}} ^{\ell}\mid O^{\ell,m}\cap L\cap
A)}{\mathbb{P}(C_{\{ \mathbf{0}\}} ^{m}\mid O^{\ell,m}\cap L\cap A)}
\mathbb{P}\bigl(C_{\{ \mathbf{0}\}} ^{m}\mid O^{\ell,m}\cap L\cap A\bigr)\\
&&\qquad = c_1^{\ell,m}\bigl(1-\mathbb{P}\bigl(C_{\{ \mathbf{0}\}} ^{\ell}\mid O^{\ell
,m}\cap L\cap A\bigr)\bigr)
\end{eqnarray*}
and therefore
%
%e11 ###
\begin{equation}\label{2star}
\mathbb{P}\bigl(C_{\{ \mathbf{0}\}} ^{\ell}\mid O^{\ell,m}\cap L\cap A\bigr)=\frac
{c_1^{\ell,m}}{c_1^{\ell,m}+1}.
\end{equation}
By similar considerations, we obtain that
%
%e12 ###
\begin{equation}\label{3star}
\mathbb{P}\bigl(C_{\{ \mathbf{0}\}} ^{\ell}\mid O^{\ell,m}\cap L\cap
A^c\bigr)=\frac{c_2^{\ell,m}}{c_2^{\ell,m}+1}
\end{equation}
and that
%
%e13 ###
\begin{equation} \label{4star}
\mathbb{P}\bigl(C_{\{ \mathbf{0}\}} ^{\ell}\mid O^{\ell,m}\cap M\bigr) =
\mathbb{P}\bigl(C_{\{ \mathbf{0}\}} ^{\ell}\mid O^{\ell,m}\cap M\cap
A^c\bigr)
= \frac{c_2^{\ell,m}}{c_2^{\ell,m}+1}.
\end{equation}

Using (\ref{2star}) and (\ref{3star}), we get that
%
%e14 ###
\begin{eqnarray} \label{5star}
\mathbb{P}\bigl(C_{\{ \mathbf{0}\}} ^{\ell}\mid O^{\ell,m}\cap L\bigr) & = & \frac
{c_1^{\ell,m}}{c_1^{\ell,m}+1}\mathbb{P}(A \mid O^{\ell,m}\cap
L)\nonumber\\
& &{} +\frac{c_2^{\ell,m}}{c_2^{\ell,m}+1}\mathbb{P}(A^c \mid O^{\ell
,m}\cap L)\nonumber\\[-8pt]\\[-8pt]
& = & \frac{c_2^{\ell,m}}{c_2^{\ell,m}+1}+\biggl(\frac{c_1^{\ell
,m}}{c_1^{\ell,m}+1}-\frac{c_2^{\ell,m}}{c_2^{\ell,m}+1}\biggr)\nonumber\\
& &\hspace*{51pt}{} \times\mathbb{P}(A \mid O^{\ell,m}\cap L).\nonumber
\end{eqnarray}
Applying (\ref{5star}) and (\ref{4star}) in (\ref{trans}) yields that,
for any $k$, we have that
%
%e15 ###
\begin{eqnarray} \label{instead}
&& 2 \bigl| \mathbb{P}\bigl(C_{\{ \mathbf{0}\}} ^{\ell}\mid L^k\bigr)-\mathbb
{P}\bigl(C_{\{ \mathbf{0}\}} ^{\ell}\mid M^k\bigr) \bigr| +  \bigl| \mathbb{P}\bigl(C_{\{
\mathbf{0}\}} ^{m}\mid L^k\bigr)-\mathbb{P}\bigl(C_{\{ \mathbf{0}\}} ^{m}\mid M^k\bigr)
\bigr| \nonumber\\[-8pt]\\[-8pt]
&&\qquad \geq \delta \biggl| \frac{c_1^{\ell,m}}{c_1^{\ell,m}+1}-\frac
{c_2^{\ell,m}}{c_2^{\ell,m}+1} \biggr|
\mathbb{P}(A \mid O^{\ell,m}\cap L^k).\nonumber
\end{eqnarray}

Since $a_{\ell}\neq1/q$ by definition, we have that $c_1^{\ell
,m}\neq c_2^{\ell,m}$.
This implies that the first two factors on the right-hand side of (\ref
{instead}) are positive constants, neither of which depends on $k$.
Now, suppose that $\mu^{\mathbb{Z}^d}_{p,q,(a_1,\ldots,a_s)}$ is
$k$-Markovian for some~$k$. In that case, the left-hand side of (\ref
{instead}) is $0$
since $\sigma^{k,\ell}_{\Lambda_k\setminus\{\mathbf{0}\}}=\sigma
^{k,m} _{\Lambda_k\setminus\{\mathbf{0}\}}$, therefore
$\mathbb{P}(A \mid O^{\ell,m}\cap L^k)=0$. This proves part 1 of Lemma
\ref{keylemma}.
Similarly,\vspace*{1pt} if $\mu^{\mathbb{Z}^d}_{p,q,(a_1,\ldots,a_s)}$ is
quasilocal, then the limit of the left-hand side of (\ref{instead}) is $0$
as $k\to\infty$, which cannot be the case if $\mathbb{P}(A \mid
O^{\ell,m}\cap L^k)$ is bounded away from $0$, uniformly in $k$.
This concludes the proof of part 2.
\end{pf}
\begin{pf*}{Proof of Theorem \ref{mainthm}, parts 1, 2\textup{(b)} and
3\textup{(b)}}
For this proof, recall the notion of
an increasing event on $\Omega_D$ (see Section \ref{rcmeasures}).
Let $d\geq2$, $q\geq1$ and $p$, $a_1,a_2,\ldots,a_s\in(0,1)$ be
arbitrary parameters of the $\operatorname{DaC}(q)$ model on $\mathbb{Z}^d$
in such a way that $S\neq S_{1/q}$ and let $\sigma^{k,\ell}$, $A$ and
$O^{\ell,m}$ be as in Lemma \ref{keylemma}.
For $k\in\{1,2,\ldots\}$, define the edge sets $E^{d,k}=\{e\in
\mathcal{E}^d\dvtx e$ is incident to $\mathbf{0}$ or to some $v\in\mathbb
{Z}^d\setminus\{\mathbf{0}\}$ with $\sigma^{k,\ell}(v)=m\}$.
For a parameter $\tilde{p}\in(0,1)$ and each $k\in\{1,2,\ldots\}$,
we define
an inhomogeneous bond percolation measure $P_{\tilde{p},k}$ on $\Omega_D$
which assigns value $0$ to all $e\in E^{d,k}$ and, independently to
each $e\in\mathcal{E}^d\setminus E^{d,k}$, value $1$ with probability
$\tilde{p}$ and $0$ with probability $1-\tilde{p}$.
It follows from Lemma \ref{condedgezd}
and Definition \ref{FKdef} that the marginal on $Y$ of the conditional
distribution $\mathbb{P}^{\mathbb{Z}^d}_{p,q,(a_1,\ldots,a_s)}$
given $O^{\ell,m}\cap C_{\mathbb{Z}^d\setminus\{\mathbf{0}\}}^{\sigma
^{k,\ell}}$
is stochastically larger than $P_{\tilde{p},k}$ with $\tilde{p}=\frac
{p}{p+(1-p)qa_{\ell}}$ if
$qa_{\ell}\geq1$, and with $\tilde{p}=p$ if $qa_{\ell}\leq1$.
Therefore, denoting the projection of $A$ on $\Omega_D$ by $A_D$ (note
that $A_D\subset\Omega_D$ is increasing), we have, for any $k$, that
%
%e16 ###
\begin{equation}\label{lowb}
\mathbb{P}^{\mathbb{Z}^d}_{p,q,(a_1,\ldots,a_s)}\bigl(A \mid O^{\ell
,m}\cap C_{\mathbb{Z}^d\setminus\{\mathbf{0}\}}^{\sigma^{k,\ell}}\bigr)\geq
P_{\tilde{p},k}(A_D).
\end{equation}
Since $\min\{\frac{p}{p+(1-p)qa_{\ell}},p\}>0$ for all $p\in(0,1)$,
we obviously have, for any fixed $k\in\{1,2,\ldots\}$, that
$P_{\tilde{p},k}(A_D)>0$.
This proves non-$k$-Markovianness of the measure $\mu^{\mathbb
{Z}^d}_{p,q,(a_1,\ldots,a_s)}$
according to (\ref{lowb}) and part 1 of Lemma \ref{keylemma}.

For the proof of part 2(b), recall
the definition of the vertices $u_1,u_2\in\mathbb{Z}^d$ and edges
$e_1,e_2,\ldots,e_{2d}\in\mathcal{E}^d$ (immediately before Lemma \ref
{keylemma}),
and the fact that $\mathcal{H}^+$ denotes the set of vertices in
$\mathbb{Z}^d$ whose first coordinate is strictly positive.
Define $\mathcal{H}^{-}$ as the set of vertices in $\mathbb{Z}^d$ whose
first coordinate is strictly negative.
Consider the events
$A_{\mathcal{H}^+}=\{\eta\in\Omega_D\dvtx\exists$ an infinite open
path in
$\eta_{\mathcal{H}^+}$ which contains $u_1\}$,
$A_{\mathcal{H}^{-}}=\{\eta\in\Omega_D\dvtx\exists$ an infinite open
path in
$\eta_{\mathcal{H}^{-}}$ which contains $u_{2}\}$,
$U=\{\eta\in\Omega_D$: there is at most one infinite open cluster in
$\eta_{\mathcal{E}^d\setminus\{e_1,e_2,\ldots,e_{2d}\}}\}$
and note that
$A_{\mathcal{H}^+}\cap A_{\mathcal{H}^{-}}\cap U\subset A_D$.

Now, assume that $a_{\ell}>1/q$. This implies that $qa_{\ell}\geq1$
and hence the free random-cluster measure $\Phi^{\mathbb
{Z}^d,0}_{p,qa_{\ell}}$ exists and is the stochastically smallest
random-cluster measure for $\mathbb{Z}^d$ with parameters $p$ and
$qa_{\ell}$ (see Section \ref{rcmeasures}).
Let us denote by $\Phi_k^{(c)}$
the measure $\Phi^{\mathbb{Z}^d,0}_{p,qa_{\ell}}$ conditioned on
the event $\{\eta\in\Omega_D\dvtx\eta_{E^{d,k}}\equiv0\}$.
Due to the aforementioned extremality of $\Phi^{\mathbb
{Z}^d,0}_{p,qa_{\ell}}$ with respect to stochastic ordering,
Lemma \ref{condedgezd} implies that the marginal on $Y$ of the measure
$\mathbb{P}^{\mathbb{Z}^d}_{p,q,(a_1,\ldots,a_s)}$
conditioned on
$O^{\ell,m}\cap C_{\mathbb{Z}^d\setminus\{\mathbf{0}\}}^{\sigma^{k,\ell
}}$ is stochastically larger than $\Phi_k^{(c)}$.
Therefore,
we have that
%
%e17 ###
\begin{eqnarray} \label{computation1}
\mathbb{P}^{\mathbb{Z}^d}_{p,q,(a_1,\ldots,a_s)}\bigl(A\mid O^{\ell,m}\cap
C_{\mathbb{Z}^d\setminus\{\mathbf{0}\}}^{\sigma^{k,\ell}}\bigr) & \geq& \Phi
_k^{(c)}(A_D) \nonumber\\[-8pt]\\[-8pt]
& \geq& \Phi_k^{(c)}(A_{\mathcal{H}^+}\cap A_{\mathcal{H}^{-}}\cap
U).\nonumber
\end{eqnarray}
Under the measure $\Phi^{\mathbb{Z}^d,0}_{p,qa_{\ell}}$, the event
$U$ has probability $1$ and the event one conditions on to obtain $\Phi
_k^{(c)}$ has positive probability, so it follows that $\Phi_k^{(c)}(U)=1$.
Hence, we have that
%
%e18 ###
\begin{eqnarray} \label{computation2}
\Phi_k^{(c)}(A_{\mathcal{H}^+}\cap A_{\mathcal{H}^{-}}\cap U) & = &
\Phi_k^{(c)}(A_{\mathcal{H}^+}\cap A_{\mathcal{H}^{-}})
\nonumber\\[-8pt]\\[-8pt]
& \geq& \Phi_k^{(c)}(A_{\mathcal{H}^+})\Phi_k^{(c)}(A_{\mathcal
{H}^{-}}),\nonumber
\end{eqnarray}
by (\ref{strongFKGfree}), since $A_{\mathcal{H}^+}$ and $A_{\mathcal
{H}^{-}}$ are increasing events.

Recalling from Section \ref{rcmeasures} that
$\tilde{E}\subset\mathcal{E}^d$ is the set of edges
that are incident to at least one vertex in $\mathbb{Z}^d\setminus
\mathcal{H}^+$, we have, again by the FKG inequality, that
%
%e19 ###
\begin{equation}\label{plus}
\Phi_k^{(c)}(A_{\mathcal{H}^+})\geq\Phi^{\mathbb{Z}^d,0}_{p,qa_{\ell
}}(A_{\mathcal{H}^+}\mid\{\eta\in\Omega_D\dvtx\eta_{\tilde{E}}\equiv
0\}).
\end{equation}
Similarly, defining the half-space $\mathcal{H}^{\geq-1}$ as the set
of vertices in $\mathbb{Z}^d$ whose first coordinate is at least $-1$
and denoting by $\tilde{E}^{\prime}$ the set of edges in $\mathcal
{E}^d$ that are incident to a vertex in $\mathcal{H}^{\geq-1}\setminus
\{u_2\}$, we have, by the FKG inequality, that
%
%e20 ###
\begin{eqnarray}\label{minus}
\Phi_k^{(c)}(A_{\mathcal{H}^-}) & \geq& \Phi^{\mathbb
{Z}^d,0}_{p,qa_{\ell}}(A_{\mathcal{H}^-}\mid\{\eta\in\Omega_D\dvtx
\eta_{\tilde{E}^{\prime}}\equiv0\}) \nonumber\\[-8pt]\\[-8pt]
& = & \frac{p}{p+(1-p)qa_{\ell}}\Phi^{\mathbb{Z}^d,0}_{p,qa_{\ell
}}(A_{\mathcal{H}^+}\mid\{\eta\in\Omega_D\dvtx\eta_{\tilde{E}}\equiv
0\}).\nonumber
\end{eqnarray}
Here, we have also used the facts that, conditioning on $\{\eta\in
\Omega_D\dvtx\eta_{\tilde{E}^{\prime}}\equiv0\}$, $A_{\mathcal{H}^-}$
can occur only if the edge between $u_2$ and $(-2,0,\ldots,0)$ is open
(which has conditional probability $\frac{p}{p+(1-p)qa_{\ell}}$ by
Definition \ref{FKdef}) and that the states of edges incident to $u_2$
are conditionally independent of the event that $(-2,0,\ldots,0)$ is
in an infinite open edge component in the corresponding half-space.
It follows from (\ref{plus}), (\ref{minus}) and the definition of
$p^{\mathcal{H}}_c(qa_{\ell},d)$ that
for all $p>p^{\mathcal{H}}_c(qa_{\ell},d)$, both
$\Phi_k^{(c)}(A_{\mathcal{H}^+})$ and $\Phi_k^{(c)}(A_{\mathcal
{H}^{-}})$ are bounded away from $0$, uniformly in $k$.
Therefore, by (\ref{computation1}) and~(\ref{computation2}),
$\mathbb{P}^{\mathbb{Z}^d}_{p,q,(a_1,\ldots,a_s)}(A\mid O^{\ell
,m}\cap C_{\mathbb{Z}^d\setminus\{\mathbf{0}\}}^{\sigma^{k,\ell}})$ is
bounded away from $0$ for such values of $p$, which implies
nonquasilocality of $\mu^{\mathbb{Z}^d}_{p,q,(a_1,\ldots,a_s)}$, by
part 2 of Lemma \ref{keylemma}.
This concludes the proof of part 2(b).

In the case where $a_{\ell}<1/q$, as remarked above,
(\ref{lowb}) holds with $\tilde{p}=p$.
On the other hand, if $p>p_c(1,d)$, then $\tilde{p}>p_c(1,d)$ and,
hence, by Lemma 8.2 in \cite{fuzzyGibbs} [whose proof is based on a
computation similar to (\ref{computation1}) and (\ref{computation2})],
we have
\[
\lim_{k\to\infty}P_{\tilde{p},k}(A_D)>0.
\]
By this, (\ref{lowb}) and part 2 of Lemma \ref{keylemma}, it follows
that $\mu^{\mathbb{Z}^d}_{p,q,(a_1,\ldots,a_s)}$ is not quasilocal,
proving part 3(b).
\end{pf*}

Our proof\vspace*{2pt} of part 2(a) of Theorem \ref{mainthm}, that is, quasilocality
of $\mu^{\mathbb{Z}^d}_{p,q,(a_1,\ldots,a_s)}$ for small $p$ when
$a_{\ell}>1/q$, will be a straightforward generalization of the proof
of part (i) of Theorem 4.4
in \cite{fuzzyGibbs}.
Although slightly more care is required when $a_{\ell}<1/q$, a similar
argument will also work in that case. Therefore, we will be able to
provide a proof below which deals with both cases simultaneously.
\begin{pf*}{Proof of Theorem \ref{mainthm}, parts 2\textup{(a)} and
3\textup{(a)}}
For the proof, recall the definitions of $\hat Y$ and $\hat\Omega$
(Definition \ref{Phatdef}), and, for a set $W\subset\mathbb{Z}^d$ and
a spin configuration $\kappa\in S^W$, recall the definition of
$K_W^{\kappa}$ (Section \ref{mainresults}) and define the analogous
event $\hat C^{\kappa}_W=\{(\xi,\eta,\hat{\eta})\in\hat{\Omega
}\dvtx\xi_W=\kappa\}$.
Fix parameters $d\geq2$, $q\geq1$ and $a_1,a_2,\ldots,a_s\in(0,1)$
of the $\operatorname{DaC}(q)$
model on $\mathbb{Z}^d$, and $p$ in such a way that
$p<p_c(qa_{\ell},d)$ if $a_{\ell}>1/q$, and
$p<\frac{p_c(1,d)qa_{\ell}}{p_c(1,d)qa_{\ell}+1-p_c(1,d)}$ if $a_{\ell}<1/q$.
Fix an arbitrary $W\subset\subset\mathbb{Z}^d$, $\kappa\in S^{W}$
and $\varepsilon>0$.
We shall show the existence of $N=N(\varepsilon,W)$ such that for all
$n\geq N$,
if $\sigma,\sigma^{\prime}\in S^{\mathbb{Z}^d\setminus W}$ are spin
configurations that agree on $\Lambda_n\setminus W$, then
%
%e21 ###
\begin{equation}\label{quasishow}
\bigl|\mu^{\mathbb{Z}^d}_{p,q,(a_1,\ldots,a_s)}(K_{W}^{\kappa} \mid
K_{\mathbb{Z}^d\setminus W}^{\sigma})-
\mu^{\mathbb{Z}^d}_{p,q,(a_1,\ldots,a_s)}(K_{W}^{\kappa} \mid
K_{\mathbb{Z}^d\setminus W}^{\sigma^{\prime}})\bigr|
\end{equation}
is less than or equal to $\varepsilon$.

In order to find such an $N$, we consider a ``dominating measure''
$\phi^{\mathrm{dom}}$ on $\Omega_D$:
we define $\phi^{\mathrm{dom}}=\Phi^{\mathbb{Z}^d,1}_{p,qa_{\ell}}$
in the case where $a_{\ell}>1/q$ and $\phi^{\mathrm{dom}}=\Phi
^{\mathbb{Z}^d,1}_{{p}/({p+(1-p)qa_{\ell}}),1}$ in the case\vspace*{2pt} where
$a_{\ell}<1/q$.
By Corollary \ref{stochlargercor},
$\phi^{\mathrm{dom}}$ is stochastically larger than\vspace*{1pt} the conditional
distribution of the modified random edge configuration $\hat{Y}$ given
any spin configuration.
Note that the parameters are chosen in such a way that
$\phi^{\mathrm{dom}}$-a.s. there exists no infinite open edge cluster [for the case
$a_{\ell}<1/q$, note that $p<\frac{p_c(1,d)qa_{\ell}}{p_c(1,d)qa_{\ell
}+1-p_c(1,d)}$ ensures that $\frac{p}{p+(1-p)qa_{\ell}}< p_c(1,d)$].
Therefore, it is possible to choose an $N$ so large that
%
%e22 ###
\begin{equation}\label{noinf}
\phi^{\mathrm{dom}}(\{\partial W\leftrightarrow\partial\Lambda_N\}
)\leq\varepsilon,
\end{equation}
where $\{\partial W\leftrightarrow\partial\Lambda_N\}=\{\eta\in
\Omega_D$: there exists a path between $\partial W$ and
$\partial\Lambda_N$ along which all edges are open in $\eta\}$.
Fix an arbitrary $n\geq N$ and let $\sigma,\sigma^{\prime} \in
S^{\mathbb{Z}^d\setminus W}$ be two arbitrary spin configurations such
that $\sigma_{\Lambda_n\setminus W}=\sigma^{\prime}_{\Lambda
_n\setminus W}$. An informal overview of the proof that $\mbox{(\ref
{quasishow})}\leq\varepsilon$ is as follows.

Let $\hat{Y}$ (resp., $\hat{Y}^{\prime}$) be the modified random edge
configuration when the spin configuration $\sigma$ (resp., $\sigma
^{\prime}$) is given.
We would like to show that $\hat{Y}$ and $\hat{Y}^{\prime}$ can be
coupled in such a way that there exists a barrier $B$ with a high
enough (at least $1-\varepsilon$) probability so that:
(a) $\hat{Y}_B=\hat{Y}^{\prime}_B\equiv0$; (b) $B$ separates
$\partial W$ and $\partial\Lambda_n$. By the definition of $\hat
{Y}$, a barrier $B$ satisfying (a) is a quasi-closed barrier in the
case where the spin configuration $\sigma$ is given. Therefore, if $B$
also satisfies (b),
then, by Corollary \ref{barrierusezd}, the spin configuration in $W$
does not depend on $\sigma_{\mathbb{Z}^d\setminus\Lambda_n}\subset
\sigma_{\mathrm{ext}(B)}$.
Clearly, the same argument holds for $\sigma^{\prime}$. Since we have
that $\sigma^{\prime}_{\Lambda_n\setminus W}=\sigma_{\Lambda
_n\setminus W}$,
we see that finding a barrier $B$ that satisfies (a) and (b) ensures
that the conditional spin distribution in $W$ is the same, given either
of $\sigma$ or $\sigma^{\prime}$.
Therefore, finding such a barrier with probability at least
$1-\varepsilon$ yields that $\mbox{(\ref{quasishow})}\leq\varepsilon$.

In order to find such a barrier, we will couple $\hat{Y}$ and $\hat
{Y}^{\prime}$ together with an auxiliary random edge configuration
$Y^{\mathrm{dom}}$ with distribution $\phi^{\mathrm{dom}}$.
We will show below how one can repeatedly use Lemma \ref{stochlarger}
to simultaneously construct $\hat{Y}$, $\hat{Y}^{\prime}$ and
$Y^{\mathrm{dom}}$ with the correct distributions in such a way that
$Y^{\mathrm{dom}}\geq\hat{Y}$ and $Y^{\mathrm{dom}}\geq\hat
{Y}^{\prime}$ hold at all stages of the construction. By the choice of
$N$ in (\ref{noinf}), we will find, with probability $1-\varepsilon$,
a barrier $B$ satisfying (b) with $Y^{\mathrm{dom}}_B\equiv0$. The
point is that since $Y^{\mathrm{dom}}\geq\hat{Y}$ and
$Y^{\mathrm{dom}}\geq\hat{Y}^{\prime}$, this implies that (a) also holds for
$B$, so we are done.
It is important to add that our construction will find an appropriate
barrier $B$ when such a barrier exists by assigning $\hat{Y}$-, $\hat
{Y}^{\prime}$- and $Y^{\mathrm{dom}}$-values only to edges in $B\cup
\mathcal{E}_{\mathrm{ext}(B)}$. Therefore,
although $\hat{Y}$ and $\hat{Y}^{\prime}$ may take different values
on such edges, the conditional spin distributions given $\sigma$ and
the explored part of $\hat Y$, respectively, given $\sigma^{\prime}$ and the
explored part of $\hat{Y}^{\prime}$ are indeed the same in $W\subset
\mathcal{V}_{\mathrm{int}(B)}$.

The formal implementation of this idea proceeds via essentially the
same coupling as used in \cite{fuzzyGibbs}, but we give it now for the
sake of completeness. We will define below a probability measure
$\mathbb{Q}$ on $\hat\Omega\times\hat\Omega\times\Omega_D$ that
is a coupling of:
\begin{longlist}
\item an $\hat\Omega$-valued random triple $(X,Y,\hat{Y})$ with
distribution $\hat{\mathbb{P}}^{\mathbb{Z}^d}_{p,q,(a_1,\ldots,a_s)}$
conditioned on $\hat C_{\mathbb{Z}^d\setminus W}^{\sigma} $;
\item an $\hat\Omega$-valued random triple $(X^{\prime
},Y^{\prime},\hat{Y}^{\prime})$ with distribution $\hat{\mathbb
{P}}^{\mathbb{Z}^d}_{p,q,(a_1,\ldots,a_s)}$ conditioned on $\hat
C_{\mathbb{Z}^d\setminus W}^{\sigma^{\prime}}$;
\item an $\Omega_D$-valued random edge configuration
$Y^{\mathrm{dom}}$ with distribution $\phi^{\mathrm{dom}}$.
\end{longlist}
It then follows from the coupling inequality (Proposition 4.2 in
\cite{GHM}) that $\mbox{(\ref{quasishow})}\leq\mathbb{Q}(X_W\neq X^{\prime}_W)$.
Hence, showing that $\mathbb{Q}(X_W=X^{\prime}_W)\geq1-\varepsilon$
would complete the proof. We define $\mathbb{Q}$ in three stages, as follows.

\begin{longlist}[II.]
\item[I.] Recall that $\mathcal{E}_n$ is the set of edges with both
endvertices in $\Lambda_n\cup\partial\Lambda_n$. It follows from
Corollary \ref{stochlargercor}
and Strassen's theorem (see Section \ref{rcmeasures}) that
the set $Q=\{\mu\dvtx\mu$ is a coupling of (i), (ii) and (iii) satisfying that
$
\mu(\hat{Y}_{\mathcal{E}^d\setminus\mathcal{E}_n}\leq
Y^{\mathrm{dom}}_{\mathcal{E}^d\setminus\mathcal{E}_n}$ and $\hat
{Y}^{\prime}_{\mathcal{E}^d\setminus\mathcal{E}_n}\leq
Y^{\mathrm{dom}}_{\mathcal{E}^d\setminus\mathcal{E}_n})=1\}
$
of probability measures on $\hat\Omega\times\hat\Omega\times
\Omega_D$ is nonempty. We will choose $\mathbb{Q}$ from this set and
will specify in stages II and III which element of $Q$ we pick.

\item[II.] Fix an arbitrary deterministic ordering of $\mathcal{E}_n$ and let
$(U_e\dvtx e\in\mathcal{E}_n)$ be a collection of independent random
variables with uniform distribution on the interval $[0,1]$.
The following algorithm will determine $\hat{Y}$, $\hat{Y}^{\prime}$
and $Y^{\mathrm{dom}}$ on a subset of $\mathcal{E}_n$, given that they
are known in $\mathcal{E}^d\setminus\mathcal{E}_n$, by drawing $\hat
{Y}$-, $\hat{Y}^{\prime}$- and $Y^{\mathrm{dom}}$-values for one edge
at a time, as follows.

\begin{enumerate}
\item Let $e\in\mathcal{E}_n$ be the first edge in the previously
fixed deterministic ordering which has not been selected in any
previous step of the algorithm and is incident to some vertex in
$\partial\Lambda_n$ or some previously selected edge $f$ with
$Y^{\mathrm{dom}}(f)=1$.
\item Let us denote by $\mathbb{P}^{(c)}$ the probability measure $\hat
{\mathbb{P}}^{\mathbb{Z}^d}_{p,q,(a_1,\ldots,a_s)}$ conditioned on
$\hat C_{\mathbb{Z}^d\setminus W}^{\sigma} $ and what we have seen
thus far of $\hat{Y}$, by $\mathbb{P}^{(c)\prime}$ the measure $\hat
{\mathbb{P}}^{\mathbb{Z}^d}_{p,q,(a_1,\ldots,a_s)}$ conditioned on
$\hat C_{\mathbb{Z}^d\setminus W}^{\sigma^{\prime}}$ and what we have
seen thus far of $\hat{Y}^{\prime}$ and by $\phi^{(c)}$ the measure
$\phi^{\mathrm{dom}}$ conditioned on what we have seen thus far of
$Y^{\mathrm{dom}}$.
We define
\[
\hat Y (e)=
\cases{
1, &\quad if $U_e<\mathbb{P}^{(c)}\bigl(\hat Y (e)=1\bigr)$,\cr
0, &\quad otherwise,}
\]
analogously,
\[
\hat Y^{\prime}(e)=
\cases{
1, &\quad if $U_e<\mathbb{P}^{(c)\prime}\bigl(\hat Y^{\prime}
(e)=1\bigr)$, \cr
0, &\quad otherwise,}
\]
and, finally,
\[
Y^{\mathrm{dom}}(e)=
\cases{
1, &\quad if $U_e<\phi^{(c)}\bigl(Y^{\mathrm{dom}}(e)=1\bigr)$, \cr
0, &\quad otherwise.}
\]
Note that if we had $Y^{\mathrm{dom}}\geq\hat{Y}$ and
$Y^{\mathrm{dom}}\geq\hat{Y}^{\prime}$ before step 2 of the algorithm (which is
$\mu$-a.s. the case for any $\mu\in Q$ before the beginning of this
algorithm), then Lemma \ref{stochlarger} implies that these
inequalities are preserved by step 2.
\item If determining $Y^{\mathrm{dom}}(e)$ in step 2 creates either
an open path in $Y^{\mathrm{dom}}$ between $\partial\Lambda_n$ and
$\partial W$ or a barrier $B$ such that $W\cup\partial W\subset
\mathcal{V}_{\mathrm{int}(B)}\subset\Lambda_n$ and $Y^{\mathrm{dom}}_B\equiv
0$, then we stop the algorithm; otherwise, we go back to step 1.
\end{enumerate}
Note that this algorithm terminates at the latest once all edges in
$\operatorname{ext}(\Delta W)$
have been selected and that it does not select any edge in $\Delta W$
or in $\mathrm{int}(\Delta W)$.

\item[III.]
If the algorithm in stage II ends by finding an open path in
$Y^{\mathrm{dom}}$ between $\partial\Lambda_n$ and $\partial W$, then we draw
the rest of $(X,Y,\hat{Y})$, $(X^{\prime},Y^{\prime},\hat
{Y}^{\prime})$ and $Y^{\mathrm{dom}}$ arbitrarily with the correct
conditional distributions, given what we have seen of them thus far.
This will possibly give that $X_W\neq X^{\prime}_W$, but that is not a
problem since, by inequality (\ref{noinf}), this case occurs with
probability at most $\varepsilon$ and, otherwise, we will always be
able to ensure that $X_W=X^{\prime}_W$.

Indeed, let as assume that the above algorithm found a barrier $B$ such
that $W\cup\partial W\subset\mathcal{V}_{\mathrm{int}(B)}\subset\Lambda_n$
and $Y^{\mathrm{dom}}_B\equiv0$. Since the inequalities
$Y^{\mathrm{dom}}\geq\hat{Y}$ and $Y^{\mathrm{dom}}\geq\hat{Y}^{\prime}$ were
retained throughout the whole algorithm (as remarked in step~2),
it follows from $Y^{\mathrm{dom}}_B\equiv0$ that $B$ is closed in
$\hat{Y}$ and $\hat{Y}^{\prime}$ as well.
Since $B$ is a barrier which is closed in $\hat Y$, it is a
quasi-closed barrier in $(X,Y)$. Therefore,
Corollary \ref{barrierusezd} implies
that the conditional distribution of $(X,Y)$ on $\mathrm{int}(B)$, given
$X_{\mathbb{Z}^d\setminus W}={\sigma}$
and what we have seen of $\hat{Y}$, is
$\mathbb{P}^{\mathrm{int}(B)}_{p,q,(a_1,\ldots,a_s)}$ conditioned on
$\{(\xi,\eta)\in\Omega^{\mathrm{int}(B)}\dvtx\xi_{\mathcal{V}_{\mathrm{int}(B)}\setminus
W}=\sigma_{\mathcal{V}_{\mathrm{int}(B)}\setminus W}\}$.
By similar considerations, the conditional distribution of $(X^{\prime
},Y^{\prime})$, given $X^{\prime}_{\mathbb{Z}^d\setminus W}={\sigma
^{\prime}}$ and what we have seen of $\hat{Y}^{\prime}$, is $\mathbb
{P}^{\mathrm{int}(B)}_{p,q,(a_1,\ldots,a_s)}$ conditioned on
$\{(\xi,\eta)\in\Omega^{\mathrm{int}(B)}\dvtx\xi_{\mathcal{V}_{\mathrm{int}(B)}\setminus
W}=\sigma^{\prime}_{\mathcal{V}_{\mathrm{int}(B)}\setminus W}\}$.
Since $\mathcal{V}_{\mathrm{int}(B)}\subset\Lambda_n$ and $\sigma_{\Lambda
_n\setminus W}=\sigma^{\prime}_{\Lambda_n\setminus W}$, we can take
$
(X_{\mathcal{V}_{\mathrm{int}(B)}},Y_{\mathcal{E}_{\mathrm{int}(B)}},\hat{Y}_{\mathcal
{E}_{\mathrm{int}(B)}})=(X^{\prime}_{\mathcal{V}_{\mathrm{int}(B)}},Y^{\prime}_{\mathcal
{E}_{\mathrm{int}(B)}},\hat{Y}^{\prime}_{\mathcal{E}_{\mathrm{int}(B)}})
$
in our coupling. This already implies that $X_W=X^{\prime}_W$ since $W
\subset\mathcal{V}_{\mathrm{int}(B)}$, so the coupling can be completed by
drawing the rest of $(X,Y,\hat{Y})$, $(X^{\prime},Y^{\prime},\hat
{Y}^{\prime})$ and $Y^{\mathrm{dom}}$ arbitrarily with the correct
conditional distributions.

These considerations yield that with a coupling $\mathbb{Q}$ of (i),
(ii) and (iii), as specified in stages I, II and III we have that
$\mathbb{Q}(X_W=X^{\prime}_W)\geq1-\varepsilon$,
which concludes the proof as noted above.\qed
\end{longlist}
\noqed\end{pf*}

The proof of Proposition \ref{almquas} is an easier application of the
concept that the existence of a (quasi-)closed barrier
``blocks the information from outside.''
Since the proof is virtually the same as the proof of Proposition 3.7
in \cite{DaC}, that is,
the analogous statement for the DaC($1$) model, we will just sketch it
for the reader's convenience.
\begin{pf*}{Proof sketch of Proposition \ref{almquas}}
Fix $W\subset\mathbb{Z}^d$ and $\sigma\in S^{\mathbb{Z}^d\setminus
W}$ such that
none of the spins in $\sigma$ percolate. By the assumption (\ref
{noperc}), this is true for almost every spin configuration.
Let $W^{\prime}\subset\mathbb{Z}^d\setminus W$ be the union of all
spin components in $\sigma_{\mathbb{Z}^d\setminus W}$ that intersect
the vertex boundary $\partial W$.
Since there is no infinite spin component in $\sigma$,
we have that $W^{\prime}$ is a finite set, hence the edge boundary
$B=\Delta(W\cup W^{\prime})$ is a closed barrier. Therefore,
it follows from Lemma \ref{barrierusezd} that the conditional
distribution of
$\mu^{\mathbb{Z}^d}_{p,q,(a_1,\ldots,a_s)}$
given $K_{\mathbb{Z}^d\setminus W}^{\sigma}$ is
$\mu^{\mathrm{int}(B)}_{p,q,(a_1,\ldots,a_s)}$ conditioned on $\{\xi\in
\Omega^{\mathrm{int}(B)}_C\dvtx\xi_{W^{\prime}}=\sigma_{W^{\prime}}\}$.

Now, recall the definition of $\partial_nW$, the $n$-neighborhood of
$W$, from Section \ref{mainresults}.
If $k$ is so large that $W^{\prime}\subset\partial_{k-1}W$, and
$\sigma^{\prime}\in S^{\mathbb{Z}^d\setminus W}$ is
such that $\sigma^{\prime} _{\partial_kW}=\sigma_{\partial_kW}$,
then it is clear, by the same argument, that
the conditional distribution of
$\mu^{\mathbb{Z}^d}_{p,q,(a_1,\ldots,a_s)}$
given $K_{\mathbb{Z}^d\setminus W}^{\sigma^{\prime}}$ is
$\mu^{\mathrm{int}(B)}_{p,q,(a_1,\ldots,a_s)}$ conditioned on $\{\xi\in
\Omega^{\mathrm{int}(B)}_C\dvtx\xi_{W^{\prime}}=\sigma^{\prime}_{W^{\prime}}\}$.
Since $\sigma_{W^{\prime}}=\sigma^{\prime}_{W^{\prime}}$, the
above conditions are the same, therefore, for any $\kappa\in S^W$, we
have that
\[
\bigl|\mu^{\mathbb{Z}^d}_{p,q,(a_1,\ldots,a_s)}(K_W^{\kappa}\mid
K_{\mathbb{Z}^d\setminus W}^{\sigma})-\mu^{\mathbb
{Z}^d}_{p,q,(a_1,\ldots,a_s)}(K_W^{\kappa}\mid K_{\mathbb
{Z}^d\setminus W}^{\sigma^{\prime}} )\bigr|=0.
\]
This proves almost sure quasilocality of $\mu^{\mathbb
{Z}^d}_{p,q,(a_1,\ldots,a_s)}$.
\end{pf*}

\section*{Acknowledgments}
I am grateful to Federico Camia and Ronald
Meester for a careful reading of an earlier version of this paper
and numerous suggestions, and to Aernout van Enter for drawing my
attention to \cite{Bodineau} and for other valuable comments.
I would also like to thank the two anonymous referees for their remarks
and corrections, which have definitely improved this paper.

% imsref loaded by lrinkeviciute, 2010-02-23 13:48:48

%
\printaddresses


\begin{thebibliography}{29}

%b1 ###
\bibitem{BCM}
\begin{barticle}[mr]
\bauthor{\bsnm{B{\'a}lint},~\bfnm{Andr{\'a}s}\binits{A.}},
  \bauthor{\bsnm{Camia},~\bfnm{Federico}\binits{F.}} \AND
  \bauthor{\bsnm{Meester},~\bfnm{Ronald}\binits{R.}}
(\byear{2009}).
\btitle{Sharp phase transition and critical behaviour in 2{D} divide and colour
  models}.
\bjournal{Stochastic Process. Appl.}
\bvolume{119}
\bpages{937--965}.
\bid{doi={10.1016/j.spa.2008.04.003}, mr={2499865}}%
\end{barticle}%
\endbibitem%

%b2 ###
\bibitem{Ising}
\begin{barticle}[vtex]
\bauthor{\bsnm{B\'alint},~\bfnm{A.}\binits{A.}},
  \bauthor{\bsnm{Camia},~\bfnm{F.}\binits{F.}} \AND
  \bauthor{\bsnm{Meester},~\bfnm{R.}\binits{R.}}
  (\byear{2010}).
\btitle{The high-temperature Ising model on the triangular lattice is a
critical Bernoulli percolation model}.
\bjournal{J. Stat. Phys.}
%Available at}
%[math.PR]}.
\bnote{To appear}.
\end{barticle}
\endbibitem

%b3 ###
\bibitem{BGN}
\begin{barticle}[mr]
\bauthor{\bsnm{Barsky},~\bfnm{David~J.}\binits{D.~J.}},
  \bauthor{\bsnm{Grimmett},~\bfnm{Geoffrey~R.}\binits{G.~R.}} \AND
  \bauthor{\bsnm{Newman},~\bfnm{Charles~M.}\binits{C.~M.}}
(\byear{1991}).
\btitle{Percolation in half-spaces: Equality of critical densities and
  continuity of the percolation probability}.
\bjournal{Probab. Theory Related Fields}
\bvolume{90}
\bpages{111--148}.
\bid{doi={10.1007/BF01321136}, mr={1124831}}
\end{barticle}
\endbibitem

%b4 ###
\bibitem{Bodineau}
\begin{barticle}[mr]
\bauthor{\bsnm{Bodineau},~\bfnm{T.}\binits{T.}}
(\byear{2005}).
\btitle{Slab percolation for the {I}sing model}.
\bjournal{Probab. Theory Related Fields}
\bvolume{132}
\bpages{83--118}.
\bid{doi={10.1007/s00440-004-0391-6}, mr={2136868}}
\end{barticle}
\endbibitem

%b5 ###
\bibitem{vEFS2}
\begin{barticle}[auto:SpringerTagBib|2009-01-14|16:51:27]
\bauthor{\bparticle{van }\bsnm{Enter},~\bfnm{A.~C.~D.}\binits{A.~C.~D.}},
  \bauthor{\bsnm{Fern\'andez},~\bfnm{R.}\binits{R.}} \AND
  \bauthor{\bsnm{Sokal},~\bfnm{A.~D.}\binits{A.~D.}}
(\byear{1991}).
\btitle{Renormalization transformations in the vicinity of first-order
  phase transitions: What can and cannot go wrong}.
\bjournal{Phys. Rev. Lett.}
\bvolume{66}
\bpages{3253--3256}.
\bid{pmid={10043740}}
\end{barticle}
\endbibitem

%b6 ###
\bibitem{vEFS}
\begin{barticle}[auto:SpringerTagBib|2009-01-14|16:51:27]
\bauthor{\bparticle{van }\bsnm{Enter},~\bfnm{A.~C.~D.}\binits{A.~C.~D.}},
  \bauthor{\bsnm{Fern\'andez},~\bfnm{R.}\binits{R.}} \AND
  \bauthor{\bsnm{Sokal},~\bfnm{A.~D.}\binits{A.~D.}}
(\byear{1993}).
\btitle{Regularity properties of position-space renormalization group
  transformations: Scope and limitations of Gibbsian theory}.
\bjournal{J. Stat. Phys.}
\bvolume{72}
\bpages{879--1167}.
\end{barticle}
\endbibitem

%b7 ###
\bibitem{vEMSS}
\begin{bincollection}[mr]
\bauthor{\bparticle{van }\bsnm{Enter},~\bfnm{A.}\binits{A.}},
  \bauthor{\bsnm{Maes},~\bfnm{C.}\binits{C.}},
  \bauthor{\bsnm{Schonmann},~\bfnm{R.~H.}\binits{R.~H.}} \AND
  \bauthor{\bsnm{Shlosman},~\bfnm{S.}\binits{S.}}
(\byear{2000}).
\btitle{The {G}riffiths singularity random field}.
In \bbooktitle{On {D}obrushin's Way. {F}rom Probability Theory to Statistical
  Physics}.
\bseries{American Mathematical Society Translations, Series 2}
\bvolume{198}
\bpages{51--58}.
\bpublisher{Amer. Math. Soc.}, \baddress{Providence, RI}.
\bid{mr={1766342}}
\end{bincollection}
\endbibitem

%b8 ###
\bibitem{FKG}
\begin{barticle}[mr]
\bauthor{\bsnm{Fortuin},~\bfnm{C.~M.}\binits{C.~M.}},
  \bauthor{\bsnm{Kasteleyn},~\bfnm{P.~W.}\binits{P.~W.}} \AND
  \bauthor{\bsnm{Ginibre},~\bfnm{J.}\binits{J.}}
(\byear{1971}).
\btitle{Correlation inequalities on some partially ordered sets}.
\bjournal{Comm. Math. Phys.}
\bvolume{22}
\bpages{89--103}.
\bid{mr={0309498}}
\end{barticle}
\endbibitem

%b9 ###
\bibitem{Garet}
\begin{barticle}[mr]
\bauthor{\bsnm{Garet},~\bfnm{Olivier}\binits{O.}}
(\byear{2001}).
\btitle{Limit theorems for the painting of graphs by clusters}.
\bjournal{ESAIM Probab. Stat.}
\bvolume{5}
\bpages{105--118}.
\bid{doi={10.1051/ps:2001104}, mr={1875666}}
\end{barticle}
\endbibitem

%b10 ###
\bibitem{Georgii}
\begin{bbook}[mr]
\bauthor{\bsnm{Georgii},~\bfnm{Hans-Otto}\binits{H.-O.}}
(\byear{1988}).
\btitle{Gibbs Measures and Phase Transitions}.
\bseries{de Gruyter Studies in Mathematics}
\bvolume{9}.
\bpublisher{de Gruyter}, \baddress{Berlin}.
\bid{mr={956646}}
\end{bbook}
\endbibitem

%b11 ###
\bibitem{GHM}
\begin{bincollection}[mr]
\bauthor{\bsnm{Georgii},~\bfnm{Hans-Otto}\binits{H.-O.}},
  \bauthor{\bsnm{H{\"a}ggstr{\"o}m},~\bfnm{Olle}\binits{O.}} \AND
  \bauthor{\bsnm{Maes},~\bfnm{Christian}\binits{C.}}
(\byear{2001}).
\btitle{The random geometry of equilibrium phases}.
In \bbooktitle{Phase Transitions and Critical Phenomena}.
\bseries{Phase Transit. Crit. Phenom.}
\bvolume{18}
\bpages{1--142}.
\bpublisher{Academic Press}, \baddress{London}.
\bid{doi={10.1016/S1062-7901(01)80008-2}, mr={2014387}}
\end{bincollection}
\endbibitem

%b12 ###
\bibitem{GrGr}
\begin{bmisc}[auto:SpringerTagBib|2009-01-14|16:51:27]
\bauthor{\bsnm{Graham},~\bfnm{B.}\binits{B.}} \AND
  \bauthor{\bsnm{Grimmett},~\bfnm{G.}\binits{G.}}
(\byear{2009}).
\bhowpublished{Sharp thresholds for the random-cluster and Ising
models. Available at}
\href{http://arxiv.org/abs/arXiv:0903.1501v1}{arXiv:0903.1501v1 [math.PR]}.
\bnote{Preprint}.
\end{bmisc}
\endbibitem

%b13 ###
\bibitem{Grimmett}
\begin{bbook}[mr]
\bauthor{\bsnm{Grimmett},~\bfnm{Geoffrey}\binits{G.}}
(\byear{1999}).
\btitle{Percolation},
\bedition{2nd} ed.
\bseries{Grundlehren der Mathematischen Wissenschaften [Fundamental Principles
  of Mathematical Sciences]}
\bvolume{321}.
\bpublisher{Springer}, \baddress{Berlin}.
\bid{mr={1707339}}
\end{bbook}
\endbibitem

%b14 ###
\bibitem{grimmett2}
\begin{bbook}[mr]
\bauthor{\bsnm{Grimmett},~\bfnm{Geoffrey}\binits{G.}}
(\byear{2006}).
\btitle{The Random-Cluster Model}.
\bseries{Grundlehren der Mathematischen Wissenschaften [Fundamental Principles
  of Mathematical Sciences]}
\bvolume{333}.
\bpublisher{Springer}, \baddress{Berlin}.
\bid{doi={10.1007/978-3-540-32891-9}, mr={2243761}}
\end{bbook}
\endbibitem

%b15 ###
\bibitem{RCrepresentations}
\begin{barticle}[mr]
\bauthor{\bsnm{H{\"a}ggstr{\"o}m},~\bfnm{O.}\binits{O.}}
(\byear{1998}).
\btitle{Random-cluster representations in the study of phase transitions}.
\bjournal{Markov Process. Related Fields}
\bvolume{4}
\bpages{275--321}.
\bid{mr={1670023}}
\end{barticle}
\endbibitem

%b16 ###
\bibitem{fuzzyposcorr}
\begin{barticle}[mr]
\bauthor{\bsnm{H{\"a}ggstr{\"o}m},~\bfnm{Olle}\binits{O.}}
(\byear{1999}).
\btitle{Positive correlations in the fuzzy {P}otts model}.
\bjournal{Ann. Appl. Probab.}
\bvolume{9}
\bpages{1149--1159}.
\bid{doi={10.1214/aoap/1029962867}, mr={1728557}}
\end{barticle}
\endbibitem

%b17 ###
\bibitem{DaC}
\begin{barticle}[mr]
\bauthor{\bsnm{H{\"a}ggstr{\"o}m},~\bfnm{Olle}\binits{O.}}
(\byear{2001}).
\btitle{Coloring percolation clusters at random}.
\bjournal{Stochastic Process. Appl.}
\bvolume{96}
\bpages{213--242}.
\bid{doi={10.1016/S0304-4149(01)00115-6}, mr={1865356}}
\end{barticle}
\endbibitem

%b18 ###
\bibitem{fuzzyGibbs}
\begin{barticle}[mr]
\bauthor{\bsnm{H{\"a}ggstr{\"o}m},~\bfnm{Olle}\binits{O.}}
(\byear{2003}).
\btitle{Is the fuzzy {P}otts model {G}ibbsian?}
\bjournal{Ann. Inst. H. Poincar\'e Probab. Statist.}
\bvolume{39}
\bpages{891--917}.
\bid{doi={10.1016/S0246-0203(03)00026-8}, mr={1997217}}
\end{barticle}
\endbibitem

%b19 ###
\bibitem{HK}
\begin{barticle}[mr]
\bauthor{\bsnm{H{\"a}ggstr{\"o}m},~\bfnm{O.}\binits{O.}} \AND
  \bauthor{\bsnm{K{\"u}lske},~\bfnm{C.}\binits{C.}}
(\byear{2004}).
\btitle{Gibbs properties of the fuzzy {P}otts model on trees and in mean
  field}.
\bjournal{Markov Process. Related Fields}
\bvolume{10}
\bpages{477--506}.
\bid{mr={2097868}}
\end{barticle}
\endbibitem

%b20 ###
\bibitem{Israel}
\begin{barticle}[mr]
\bauthor{\bsnm{Israel},~\bfnm{R.~B.}\binits{R.~B.}}
(\byear{2004}).
\btitle{Some generic results in mathematical physics}.
\bjournal{Markov Process. Related Fields}
\bvolume{10}
\bpages{517--521}.
\bid{mr={2097870}}
\end{barticle}
\endbibitem

%b21 ###
\bibitem{KW}
\begin{barticle}[mr]
\bauthor{\bsnm{Kahn},~\bfnm{Jeff}\binits{J.}} \AND
  \bauthor{\bsnm{Weininger},~\bfnm{Nicholas}\binits{N.}}
(\byear{2007}).
\btitle{Positive association in the fractional fuzzy {P}otts model}.
\bjournal{Ann. Probab.}
\bvolume{35}
\bpages{2038--2043}.
\bid{doi={10.1214/009117907000000042}, mr={2353381}}
\end{barticle}
\endbibitem

%b22 ###
\bibitem{Kulske1}
\begin{barticle}[mr]
\bauthor{\bsnm{K{\"u}lske},~\bfnm{C.}\binits{C.}}
(\byear{1999}).
\btitle{({N}on-){G}ibbsianness and phase transitions in random lattice spin
  models}.
\bjournal{Markov Process. Related Fields}
\bvolume{5}
\bpages{357--383}.
\bid{mr={1734240}}
\end{barticle}
\endbibitem

%b23 ###
\bibitem{Kulske2}
\begin{barticle}[mr]
\bauthor{\bsnm{K{\"u}lske},~\bfnm{Christof}\binits{C.}}
(\byear{2001}).
\btitle{Weakly {G}ibbsian representations for joint measures of quenched
  lattice spin models}.
\bjournal{Probab. Theory Related Fields}
\bvolume{119}
\bpages{1--30}.
\bid{doi={10.1007/PL00012737}, mr={1813038}}
\end{barticle}
\endbibitem

%b24 ###
\bibitem{KO}
\begin{barticle}[mr]
\bauthor{\bsnm{K{\"u}lske},~\bfnm{Christof}\binits{C.}} \AND
  \bauthor{\bsnm{Opoku},~\bfnm{Alex~A.}\binits{A.~A.}}
(\byear{2008}).
\btitle{The posterior metric and the goodness of {G}ibbsianness for transforms
  of {G}ibbs measures}.
\bjournal{Electron. J. Probab.}
\bvolume{13}
\bpages{1307--1344}.
\bid{mr={2438808}}
\end{barticle}
\endbibitem

%b25 ###
\bibitem{MVV}
\begin{barticle}[mr]
\bauthor{\bsnm{Maes},~\bfnm{Christian}\binits{C.}} \AND
  \bauthor{\bsnm{Vande~Velde},~\bfnm{Koen}\binits{K.}}
(\byear{1995}).
\btitle{The fuzzy {P}otts model}.
\bjournal{J. Phys. A}
\bvolume{28}
\bpages{4261--4270}.
\bid{mr={1351929}}
\end{barticle}
\endbibitem

%b26 ###
\bibitem{NSch}
\begin{barticle}[mr]
\bauthor{\bsnm{Newman},~\bfnm{C.~M.}\binits{C.~M.}} \AND
  \bauthor{\bsnm{Schulman},~\bfnm{L.~S.}\binits{L.~S.}}
(\byear{1981}).
\btitle{Infinite clusters in percolation models}.
\bjournal{J.~Stat. Phys.}
\bvolume{26}
\bpages{613--628}.
\bid{mr={648202}}
\end{barticle}
\endbibitem

%b27 ###
\bibitem{Pisztora}
\begin{barticle}[mr]
\bauthor{\bsnm{Pisztora},~\bfnm{\'{A}goston}\binits{\'{A}.}}
(\byear{1996}).
\btitle{Surface order large deviations for {I}sing, {P}otts and percolation
  models}.
\bjournal{Probab. Theory Related Fields}
\bvolume{104}
\bpages{427--466}.
\bid{doi={10.1007/BF01198161}, mr={1384040}}
\end{barticle}
\endbibitem

%b28 ###
\bibitem{Strassen}
\begin{barticle}[mr]
\bauthor{\bsnm{Strassen},~\bfnm{V.}\binits{V.}}
(\byear{1965}).
\btitle{The existence of probability measures with given marginals}.
\bjournal{Ann. Math. Statist.}
\bvolume{36}
\bpages{423--439}.
\bid{mr={0177430}}
\end{barticle}
\endbibitem

%b29 ###
\bibitem{Wouts}
\begin{barticle}[mr]
\bauthor{\bsnm{Wouts},~\bfnm{Marc}\binits{M.}}
(\byear{2008}).
\btitle{A coarse graining for the {F}ortuin--{K}asteleyn measure in random
  media}.
\bjournal{Stochastic Process. Appl.}
\bvolume{118}
\bpages{1929--1972}.
\bid{doi={10.1016/j.spa.2007.11.009}, mr={2462281}}
\end{barticle}
\endbibitem

\end{thebibliography}
\end{document}